\documentclass[11pt]{article}
\usepackage{graphicx,amssymb,mathrsfs,amsmath,color,pifont}
\newfont{\bb}{msbm10}

\newcommand{\tr}{^{\sf T}}
\newcommand{\m}[1]{{\bf{#1}}}
\newcommand{\C}[1]{{\cal {#1}}}

\newtheorem{remark}{Remark}[section]

\newtheorem{theorem}{Theorem}[section]

\newtheorem{lemma}{Lemma}[section]
\newtheorem{corollary}{Corollary}[section]

\begin{document}
\cleardoublepage \pagestyle{plain}
\bibliographystyle{plain}

\title{ A New Insight on  Augmented Lagrangian Method with Applications in Machine Learning}

\author{
Jianchao Bai
  \footnote{Research \& Development Institute of Northwestern Polytechnical University in Shenzhen, Shenzhen 518057, China;  School of Mathematics and Statistics,  Northwestern Polytechnical
        University, Xi'an  710129,      China. (\tt jianchaobai@nwpu.edu.cn).}
\quad
Linyuan Jia
 \footnote{School of Power and Energy,  Northwestern Polytechnical
        University, Xi'an  710129,      China. (\tt jialinyuan@nwpu.edu.cn).}
\quad
Zheng Peng
\footnote{ Corresponding author. School of Mathematics and Computational Science, Xiangtan University, Xiangtan 411105, China. (\tt pzheng@xtu.edu.cn).}
     }
\date{ }
\maketitle

\vskip 3mm\noindent {\small\bf Abstract.}
By exploiting  double-penalty  terms for the primal subproblem,
we develop a novel relaxed augmented Lagrangian method for solving a family
of convex optimization  problems  subject to equality or inequality
constraints. The  method is then extended to solve a general  multi-block separable
convex optimization problem, and two related primal-dual hybrid gradient algorithms
are also discussed. Convergence results about the sublinear and linear convergence
rates are established by variational characterizations  for  both the saddle-point
of the problem and the first-order optimality conditions of involved subproblems. A
large number of experiments on testing the linear  support vector machine problem and
the robust principal component analysis problem arising from machine learning indicate
that our proposed algorithms perform much better than several state-of-the-art algorithms.


\bigskip

\section{ \large Introduction}
An    interesting work, that is  the  Balanced Augmented Lagrangian Method (abbreviated by  B-ALM)    proposed by He-Yuan \cite{Heyuan21}, aims to solve the following
  convex optimization problem subject to    equality or inequality constraints:
\begin{equation} \label{Sec1-Prob1}
\min   \big\{\theta(\m{x}) |~ A\m{x}=b~ (\operatorname{or} \geq b),~ \m{x} \in \C{X} \big\},
\end{equation}
where   $\theta:\C{R}^n \rightarrow
\C{R}$ is a   closed  proper convex function; $\C{X} \subseteq \C{R}^{n} $ is a closed  convex   set; $A\in\C{R}^{m\times n}$ and $b\in\C{R}^m$ are given. Hereafter,
    the symbols    $\C{R}, \C{R}^n(\C{R}^n_+)$ and $ \C{R}^{m\times n}$
denote  the sets    of  real numbers,     $n$ dimensional real (nonnegative) column vectors, and     $m\times n$   real matrices, respectively. The bold $\m{I}$  denotes  the identity matrix and  $\m{0}$    stands for
  zero matrix/vector  with proper dimensions.  $Q\succ\m{0}$ means  $Q$  is a symmetric positive definite   matrix{, and $\nabla f(x)$ denotes the gradient of  differentiable function $f$ at $x$.}  We   use $\|\cdot\|$ and $\langle\cdot,\cdot\rangle$  to denote the standard Euclidean norm and inner product, respectively.    Given  $H\succ\m{0}$, we define  $\|\m{w}\|_H =\sqrt{\langle\m{w}, H \m{w}\rangle}$.  Throughout this paper, the  solution set of   the problem (\ref{Sec1-Prob1}) is  assumed to be nonempty.

A fundamental tool to solve the problem (\ref{Sec1-Prob1}) is the  Augmented Lagrangian Method (ALM, \cite{Hest69, Powell69}) by exploring the following two steps:
\[\left\{\begin{array}{l}
\m{x}^{k+1}=\arg \min\limits_{\m{x} \in \C{X}}  L(\m{x}, \lambda^k)+ \frac{r}{2} \|A\m{x}-b\|^2, \\
\lambda^{k+1}= \lambda^k-r(A\m{x}^{k+1}-b),
\end{array}\right.
\]
where $r>0$ denotes penalty parameter  for the violation of the linear constraints   and $
 L(\m{x}, \lambda)=\theta(\m{x})-\langle \lambda, A\m{x}-b\rangle
 $
 denotes the corresponding Lagrangian function. With simple algebra,   the core subproblem of ALM amounts to \begin{equation}\label{sub-ALM}
 \m{x}^{k+1}=\arg\min\limits_{\m{x}\in\C{X}} \theta(\m{x})-
 \big\langle \m{x}, A\tr \lambda^k +r A\tr b\big\rangle
 +\frac{1}{2}\big\|x\big\|^2_{r A\tr A}.
 \end{equation}  It  is often complicated   and  lacks an efficient solution unless inexact approximation techniques are employed.
As described in \cite{Heyuan21}, the    B-ALM reads
\[\left\{\begin{array}{l}
\m{x}^{k+1}=\arg \min\limits_{\m{x} \in \C{X}} L(\m{x}, \lambda^k)+\frac{r}{2}\left\|\m{x}-\m{x}^{k}\right\|^{2} , \\
\lambda^{k+1}=\arg\max L([2 \m{x}^{k+1}-\m{x}^{k}], \lambda)-\frac{1}{2}\left\|\lambda-\lambda^{k}\right\|^2_{\frac{1}{r} A A\tr+\delta \m{I}},
\end{array}\right.
\]
whose convergence  depends on the positive definiteness of the  block matrix
 $
\left[\begin{array}{cccccc}
  r\m{I}  &  A\tr  \\
   A   &      A A\tr/r +\delta \m{I}
\end{array}\right]
$
for any $ r,\delta>0$.
One may use a general form $Q+ A A\tr/r $  for any $Q\succ\m{0}$ to replace the above lower-right block  so as to guarantee the convergence.
 The major merit of  B-ALM is that it
 greatly weakens  the convergence conditions of some ALM and   related first-order splitting algorithms \cite{BLW20,chpo11,GHYu14,hy13,hmy20,Ma18,ZhuW15}. Namely, the     parameter $r$ does not depend on $\rho(A\tr A)$, where $\rho(\cdot)$ represents the spectrum radius of a matrix. Another merit of B-ALM is that it simplifies the   solving difficulty  of the subproblems to a relatively easier proximal estimation, which  may have a closed-form solution in many practical applications.
 However, it will    need  an inner solver to tackle the dual subproblem or   the     Cholesky factorization to deal with an equivalent linear equation of the   dual subproblem.

Motivated by these observations,    we will   develop and analyze a new double-penalty ALM with a relaxation step  (abbreviated by  \textbf{P-rALM})  to solve the problem (\ref{Sec1-Prob1}). The method could reduce  the difficulty of updating the dual variable  while still maintaining  the weak  convergence condition   and having  a  fast convergence behavior.
For  the sake of conciseness, we first present  the framework of P-rALM   as  the following.
\begin{flushleft}
\centering\fbox{
	\parbox{0.96\textwidth}{
\[
   \begin{array}{lll}
 \hspace*{-.1in}\textrm{\bf  Initialize } (\m{x}^0,\lambda^0), \textrm{ and choose }
 \gamma\in(0,2),~ r>0,~ Q\succ\m{0}; \\
 \hspace*{-.1in} \textrm{\bf While } \textrm{the stopping criterion is not satisfied, }  \textrm{\bf do}\\
\hspace*{.1in}      \widetilde{\m{x}}^k=\arg \min\limits_{\m{x} \in \C{X}} \theta(\m{x})- \big\langle \lambda^{k}, A\m{x}-b \big\rangle +\frac{r}{2}\left\|A(\m{x}-\m{x}^{k})\right\|^2 + \frac{1}{2}\left\|\m{x}-\m{x}^{k}\right\|^2_{Q};\\
\hspace*{.1in}
\widetilde{\lambda}^k= \C{P}_{\Lambda}\Big(\lambda^{k} - r \left[A(2\widetilde{\m{x}}^k-\m{x}^k)-b\right]\Big);\\
\hspace*{.1in}  \left(\begin{array}{c}
 \m{x}^{k+1}  \\   \lambda^{k+1}
\end{array}\right) =  \left(\begin{array}{c}
 \m{x}^k  \\   \lambda^k
\end{array}\right) + \gamma\left(\begin{array}{c}
\widetilde{\m{x}}^k- \m{x}^k \\  \widetilde{\lambda}^k- \lambda^k
\end{array}\right);\\
 \hspace*{-.1in} \textrm{\bf End while}
\end{array}
\]}}\end{flushleft}
 In the above framework, $\C{P}_{\Lambda}(\cdot)$ denotes the projection operator onto the set \[\Lambda=\left\{\begin{array}{lll}
  \C{R}^m, & \textrm{ if } A\m{x}=b, \\
 \C{R}^m_+, & \textrm{ if } A\m{x}\geq b.
\end{array}\right.\]
For  the problem (\ref{Sec1-Prob1}) subject to   $A\m{x}\geq b$, the   dual update  reduces to $\widetilde{\lambda}^k=\max\big(\lambda^{k} - r \big[A(2\widetilde{\m{x}}^k-\m{x}^k)-b\big],\m{0}\big)$.
      Main features of  this  P-rALM are  summarized as   four aspects:
\begin{itemize}

\item
 {Compared to} the classical ALM { and B-ALM, the  $\m{x}$-subproblem of P-rALM utilizes   different  quadratic penalty terms  $\frac{r}{2}\left\|A(\m{x}-\m{x}^{k})\right\|^2$ and $\frac{1}{2}\left\|\m{x}-\m{x}^{k}\right\|^2_{Q}$. Equivalently,}
  \begin{equation}\label{sub-PALM}
 { \widetilde{\m{x}}^k=\arg \min\limits_{\m{x} \in \C{X}}  \theta(\m{x})-
 \big\langle \m{x}, A\tr\lambda^{k}+(Q+r A\tr A)\m{x}^{k} \big\rangle +\frac{1}{2}
 \big\|\m{x}\big\|^2_{Q+r A\tr A}},
 \end{equation}
 which is   different from (\ref{sub-ALM}) since   the iterate $ \widetilde{\m{x}}^k$   does not depend on the   data $b$ but the previous iterate $\m{x}^{k}$. { Although it seems much complicated,
  taking $Q= \tau \m{I} - r A\tr A$ with $\tau> r \|A\tr A\|$ could convert (\ref{sub-PALM})   to the following    proximity operator
\[
\textrm{\bf prox}_{\theta,\tau}(x)=\arg\min\limits_{\m{x} \in \C{X}} \theta(\m{x}) +\frac{\tau}{2}\Big \|\m{x} - \m{x}^k-\frac{1}{\tau}A\tr \lambda^k\Big\|^2,
\]
which   has a  unique global solution   and further allows  a closed-form solution when $\C{X}$ is simple.}
{For this case, by taking  $\gamma=1,\tau=\frac{r}{\alpha}$ for some $\alpha>0$,   our P-rALM  reduces to   \cite[the scheme (38)]{LXZh17}, which indicates that the new    algorithm is more general than some in the literature.} If $\textrm{\bf prox}_{\theta,\tau}(x)$ is  unavailable  but $\theta$ is smooth, then one could exploit linearization techniques  or   inner solvers (e.g. the conjugate gradient method)  to solve the $\m{x}$-subproblem inexactly, or  use the  formula provided by \cite{Osher23}    for
accurately approximating the
proximal operator.
\item
{ The dual update $\widetilde{\lambda}^k$} is the same as that in \cite{chpo11} but   { is much  easier than  that of  B-ALM}. It combines  the  information of both the current iterate $\widetilde{\m{x}}^k$ and the extrapolation iterate $\widetilde{\m{x}}^k-\m{x}^k$. Moreover, after the primal-dual updates,   a relaxation step is employed to accelerate the convergence of the algorithm  from  theoretical  and numerical interests.

\item
As previously stated,  compared to some existing splitting algorithms, the global convergence of  P-rALM  will no longer depend on $\rho(A\tr A)$, although   P-rALM with the choice $Q= \tau \m{I} - r A\tr A$ reduces to   \cite[Algorithm  1]{JWZC20} with   $\eta=0$ involved (in fact, the convergence of this algorithm  depends on $\rho(A\tr A)$). {The main contribution of both B-ALM and ours is that the algorithmic parameters
are more flexible than the previous.} Besides, the double-penalty terms in P-rALM can also be extended to the general Bregman distance   while still ensuring  the convergence   under a proper assumption.
\item
We provide two elegant results   in  Theorem \ref{b-10-1} and Corollary  \ref{b-15}, that is,  the primal residual and the objective gap converge in a sublinear convergence rate.  Motivated by the structure of the involved matrix  $H$, we also discuss a generalization of P-rALM and two new types of  Primal-Dual Hybrid Gradient algorithm (PDHG) for solving   the multiple block separable convex optimization and
   the saddle-point problem, respectively. The connection between  P-rALM and its variants is that  all of them can be analyzed by   variational characterization  with similar structured matrices to   $H$ in the sense of primal-dual and dual-primal frameworks, respectively.   The   linear convergence rate of P-rALM is indicated under   similar assumptions as the   analysis for the new PDHG in the appendix. Performance of   our     P-rALM  and its two-block extensions are   verified by testing two popular examples in machine learning and by comparing with several  well-established algorithms in the literature.
\end{itemize}

The paper is organized as follows. In Section \ref{Sec2}, we analyze the global convergence and sublinear convergence rate of P-rALM. A linearized P-rALM is   discussed   when the objective function is smooth. Section \ref{Sec3} extends  the proposed P-rALM to solve a multi-block separable convex programming    and also provides a  general dual-primal version. Section \ref{Sec4} investigates the   performance of the proposed algorithm and its  extensions.   In the appendix, we further discuss  the convergence complexity of two  related  PDHG algorithms  based on the construction of P-rALM  for solving a family of convex-concave saddle-point problems.

\section{Convergence analysis of P-rALM}\label{Sec2}
\subsection{Variational characterization}
Let us first recall the following   fundamental  lemma  \cite{He2016} which will be used to characterize the saddle-point of (\ref{Sec1-Prob1}) and the iterates of P-rALM.
\begin{lemma} \label{opt-1}
Let   $\Omega\subseteq \C{R}^n$  be a   closed convex
set,    $f(x)$ and
$h(x)$ be  convex functions. If $h$ is differentiable on an open set which contains $\Omega$, and    the solution set of the minimization  problem
$\min \left\{f(x)+h(x) \mid x\in\Omega\right\}$ is nonempty. Then, we have
\[
x^* \in \arg\min \big\{f(x)+h(x) \mid x\in\Omega\big\}\]
if and only if
\[ \ x^*\in \Omega,~ \  f(x)-f(x^*)
+\big\langle x-x^*, \nabla h(x^*)\big\rangle\geq 0,~ \forall x\in\Omega.
\]
\end{lemma}

 From the perspective of optimization,   a    point
$
\m{w}^*=(\m{x}^*;\lambda^*)\in \C{M}:= \C{X}\times \Lambda
$
 is called the saddle-point of (\ref{Sec1-Prob1})
if
\[
 L\left(\m{x}^*, \lambda\right)\leq L\left(\m{x}^*, \lambda^*\right)\leq L\left(\m{x}, \lambda^*\right),\quad\forall \lambda\in \Lambda,\m{x}\in\C{X},
\]
  which, by Lemma \ref{opt-1}, is explicitly  written  as
\[
   \left \{\begin{array}{lllll}
\m{x}^*\in\C{X}, &\theta(\m{x})- \theta(\m{x}^*) &+ &\langle\m{x}-\m{x}^*,  -A\tr \lambda^*\rangle\geq  0,  &\forall \m{x} \in\C{X},\\
\lambda^*\in \Lambda, & &&\langle\lambda-\lambda^*,  A\m{x}^*- b\rangle \geq   0, & \forall \lambda \in \Lambda.
\end{array}\right.
\]
These inequalities  can be  further expressed as
the following more compact  form
\begin{equation}\label{Sec3-1}
 {\bf\textrm{VI}(\theta,\C{J},\C{M})}: ~~  \theta(\m{x})- \theta(\m{x}^*) + \big\langle \m{w}-\m{w}^*, \C{J}(\m{w}^*)\big\rangle \geq  0, \ \forall \m{w}\in \C{M},
\end{equation}
where
\[
\m{w}=\left(\begin{array}{c}
 \m{x}  \\   \lambda
\end{array}\right) \quad \textrm{and} \quad
\C{J}(\m{w})=\left(\begin{array}{c}
 -A\tr \lambda\\    A\m{x}-b
\end{array}\right).
\]
Similar characterization     can be found in  e.g.  \cite{BJXZh18,HeY18}. An equivalent expression of (\ref{Sec3-1}) is
\begin{equation}\label{b-1}
  \theta(\m{x})- \theta(\m{x}^*) +\big\langle \m{w}-\m{w}^*, \C{J}(\m{w})\big\rangle \geq  0, \ \forall \m{w}\in \C{M},
\end{equation}
since   it holds by the monotonicity of   $\C{J}(\m{w})$ that
\begin{equation}\label{b-2}
  \big\langle \m{w}-\bar{\m{w}}, \C{J}(\m{w})-\C{J}(\bar{\m{w}})\big\rangle = 0, \ \forall \m{w}, \bar{\m{w}}\in \C{M}.
\end{equation}

The solution set of (\ref{Sec3-1})  is    nonempty by the previous assumption that the  solution set of   (\ref{Sec1-Prob1}) is   nonempty,    and it  can be characterized as  (see \cite[ Theorem 2.1]{HeYuanSIAM12})
  \[
    \C{M}^*=\bigcap\limits_{\m{w}\in \C{M}}\Big\{\hat{\m{w}}\mid
   \theta(\m{x})- \theta(\hat{\m{x}}) +\big\langle \m{w}-\hat{\m{w}}, \C{J}(\m{w})\big\rangle \geq  0 \Big\}.
   \]
   Obviously,  $\m{w}^*\in \C{M}^*$ satisfies (\ref{Sec3-1}) if and only if it is a primal-dual solution of   (\ref{Sec1-Prob1}).  Next, we characterize    the sequence generated by P-rALM as a mixed variational
inequality  with the aid of the auxiliary notation $\widetilde{\m{w}}^k=(\widetilde{\m{x}}^k;\widetilde{\lambda}^k)$.

\begin{lemma} \label{x1-optimal}
The sequence  $\{ \m{w}^k \}$   generated by P-rALM satisfies
\begin{equation}\label{b-3}
\widetilde{\m{w}}^k\in \C{M},~  \theta(\m{x})- \theta(\widetilde{\m{x}}^k) +\big\langle \m{w}-\widetilde{\m{w}}^k, \C{J}(\m{w})\big\rangle \geq   \big\langle\m{w}-\widetilde{\m{w}}^k, H(\m{w}^k-\widetilde{\m{w}}^k) \big\rangle
\end{equation}
for any $\m{w}\in \C{M}$,  where
\begin{equation}\label{H}
\widetilde{\m{w}}^k=\left(\begin{array}{c}
 \widetilde{\m{x}}^k  \\   \widetilde{\lambda}^k
\end{array}\right)  \quad \textrm{and} \quad H=\left[\begin{array}{ccccccc}
  r A\tr A+Q &&  A\tr  \\
   A   &&     \frac{1}{r} \mathbf{I}
\end{array}\right]
\end{equation}
is symmetric positive definite for any $r>0$ and $Q\succ \m{0}$.
\end{lemma}
 Proof.
According to Lemma \ref{opt-1}, we have from   the first-order optimality condition  of the $\m{x}$-subproblem  that
\[
 \widetilde{\m{x}}^k\in\C{X},~ \theta(\m{x})- \theta(\widetilde{\m{x}}^k)  +  \big\langle\m{x}-\widetilde{\m{x}}^k, -A\tr \lambda^k+(rA\tr A+Q) (\widetilde{\m{x}}^k-\m{x}^k)\big\rangle\geq  0,   ~ \forall \m{x} \in\C{X},\]
 equivalently,
\begin{eqnarray}\label{b-5}
&&\quad \theta(\m{x})- \theta(\widetilde{\m{x}}^k)  +  \big\langle\m{x}-\widetilde{\m{x}}^k, -A\tr \widetilde{\lambda}^k\big\rangle\nonumber\\
&&\geq    \Big\langle\m{x}-\widetilde{\m{x}}^k,  (rA\tr A+Q) (\m{x}^k-\widetilde{\m{x}}^k)+ A\tr (\lambda^k-\widetilde{\lambda}^k)  \Big\rangle, ~ \forall \m{x} \in\C{X}.
 \end{eqnarray}
 Besides, it follows from the dual update  that $\widetilde{\lambda}^k\in \Lambda$ and
 \begin{equation}\label{b-6}
  \Big\langle\lambda-\widetilde{\lambda}^k,  A\widetilde{\m{x}}^k-b  \Big\rangle{\geq}  \Big\langle\lambda-\widetilde{\lambda}^k,   A(\m{x}^k-\widetilde{\m{x}}^k) +\frac{1}{r}(\lambda^k-\widetilde{\lambda}^k)\Big\rangle,~\forall \lambda\in \Lambda.
\end{equation}
Combine the   inequalities in (\ref{b-5})-(\ref{b-6}) together with the structure of $H$   and the property in (\ref{b-2}) to ensure the result in    (\ref{b-3}).

Observing that the symmetric matrix $H$ has the following decomposition:
\[
H=\left[\begin{array}{ccccccc}
  r A\tr A  &&  A\tr  \\
   A   &&     \frac{1}{r} \mathbf{I}
\end{array}\right]+\left[\begin{array}{ccccccc}
   Q &&  \m{0}  \\
   \m{0}   &&    \m{0}
\end{array}\right]=  \left(\begin{array}{ccccccc}
  \sqrt{r} A\tr    \\
     \frac{1}{\sqrt{r}} \mathbf{I}
\end{array}\right)\left(\sqrt{r} A, \frac{1}{\sqrt{r}} \mathbf{I}\right) +\left[\begin{array}{ccccccc}
   Q &&  \m{0}  \\
   \m{0}   &&    \m{0}
\end{array}\right].
\]
For any $\m{w}=(\m{x};\lambda)\neq \m{0}$  we have
\[
\m{w}\tr H\m{w} = \Big\|\sqrt{r} A \m{x}+ \frac{1}{\sqrt{r}} \lambda\Big\|^2 +\big\| \m{x} \big\|^2_Q>0,
\]
and therefore   $H$ is  a positive definite matrix.
 $\hfill\blacksquare$

\subsection{Convergence and  convergence rate}

The following theorem shows that the sequence  $\{ \m{w}^k \}$   generated by P-rALM  is contractive under the $H$-weighted norm and thus  converges  to the solution point of $\textrm{VI}(\theta,\C{J},\C{M})$.

\begin{theorem}  \label{Sec31-3}
For any   $\gamma\in (0,2)$, the sequence  $\{ \m{w}^k \}$   generated by P-rALM satisfies
 \begin{equation}\label{contract}
   \big \|\m{w}^{k+1}-\m{w}^*\big\|^2_H\leq\big\|\m{w}^k-\m{w}^*\big\|^2_H -\frac{2-\gamma}{\gamma}
     \big\|\m{w}^k-\m{w}^{k+1}\big\|^2_H, \ \forall \m{w}^*\in \C{M}^*.
\end{equation}
Moreover, {there} exists a point $\m{w}^\infty\in \C{M}^*$ such that  $\lim\limits_{k\rightarrow\infty} \m{w}^{k}=\m{w}^\infty.$
\end{theorem}
 Proof.
Setting $\m{w}=\m{w}^*$ in  (\ref{b-3}) together with (\ref{Sec3-1})  is to achieve
\[\big\langle\widetilde{\m{w}}^k-\m{w}^*, H(\m{w}^k-\widetilde{\m{w}}^k) \big\rangle  \geq  \theta(\widetilde{\m{x}}^k)-\theta(\m{x}^*)  +\big\langle \widetilde{\m{w}}^k-\m{w}^*, \C{J}(\m{w}^*)\big\rangle\geq 0.
\]
Combining the above property with the update
\begin{equation}\label{rela}
\m{w}^{k+1}=\m{w}^k+\gamma(\widetilde{\m{w}}^k-\m{w}^k),
\end{equation}
 we have
\begin{eqnarray}
&&\big\|\m{w}^{k}-\m{w}^*\big\|^2_H-
\big\|\m{w}^{k+1}-\m{w}^*\big\|^2_H\nonumber \\&=&
\big\|\m{w}^{k}-\m{w}^*\big\|^2_H
- \big\|\m{w}^k-\m{w}^*+\m{w}^{k+1}-\m{w}^k\big\|^2_H  \nonumber \\
& =&   2\gamma \big\langle \m{w}^k-\m{w}^*, H(\m{w}^k-\widetilde{\m{w}}^k)\big\rangle-\gamma^2
\big\|\widetilde{\m{w}}^k-\m{w}^k\big\|^2_H  \nonumber \\
&=&2\gamma \big\langle \m{w}^k-\widetilde{\m{w}}^k+\widetilde{\m{w}}^k-\m{w}^*, H(\m{w}^k-\widetilde{\m{w}}^k)\big\rangle
-\gamma^2\big\|\m{w}^k-\widetilde{\m{w}}^k\big\|^2_H  \nonumber\\
&=& \gamma(2-\gamma)\big\|\m{w}^k-\widetilde{\m{w}}^k\big\|^2_H+ 2\gamma
\big\langle
\widetilde{\m{w}}^k-\m{w}^*, H(\m{w}^k-\widetilde{\m{w}}^k)\big\rangle\nonumber\\
&\geq& \gamma(2-\gamma)\big\|\m{w}^k-\widetilde{\m{w}}^k\big\|^2_H
= \frac{2-\gamma}{\gamma}\big\|\m{w}^k-\m{w}^{k+1}\big\|^2_H.\nonumber
\end{eqnarray}
Then,  rearrange this inequality  to obtain   the inequality (\ref{contract}).

Note that the inequality (\ref{contract}) implies that the  sequence  $\{ \m{w}^k \}$  is bounded and
 $
\lim\limits_{k\rightarrow\infty} \big\|\m{w}^k-\m{w}^{k+1}\big\|^2_H=0.$ Namely, $\lim\limits_{k\rightarrow\infty}  (\m{w}^k-\m{w}^{k+1}) =\m{0},
$
which by the relation (\ref{rela}) also shows
\begin{equation}\label{b-9}
\lim\limits_{k\rightarrow\infty}(\widetilde{\m{w}}^k-\m{w}^k)=\m{0}.
\end{equation}
Let $\m{w}^\infty$ be any accumulation point of $\{\widetilde{\m{w}}^k \}$.  Then, (taking a subsequence of $\{\widetilde{\m{w}}^k \}$ if
necessary) it follows from (\ref{b-3}) and (\ref{b-9}) that
\[ \theta(\m{x})- \theta(\m{x}^\infty) +\left\langle \m{w}-\m{w}^\infty, \C{J}(\m{w}^\infty)\right\rangle  \geq 0, \ \forall \m{w}\in \C{M}.
     \]
     This indicates $\m{w}^\infty\in \C{M}^*$ compared to (\ref{Sec3-1}). So, by (\ref{contract}) again  we have
    \[
      \big\|\m{w}^k-\m{w}^\infty \big\|^2_H\leq\big\|\m{w}^j-\m{w}^\infty \big\|^2_H  ~~~~ \textrm{ for~all } k\geq j.
      \]
    Finally, it follows from  $\m{w}^\infty$ being an accumulation point that $\lim\limits_{k\rightarrow\infty} \m{w}^{k}=\m{w}^\infty.$
   $\hfill\blacksquare$

Before   establishing the  sublinear  convergence
rate of   P-rALM for the following   average iterates (it was first  used in \cite{BH20} to accelerate the convergence of a stochastic method)
\begin{equation}\label{erg-iterate}
\m{w}_T :=\frac{1}{T+1}\sum_{k=\kappa}^{T+\kappa}\widetilde{\m{w}}^k\quad \textrm{and} \quad
\m{x}_T :=\frac{1}{T+1}\sum_{k=\kappa}^{T+\kappa}\widetilde{\m{x}}^k,\quad {\forall \kappa\geq 0, T>0,}
\end{equation}
we    analyze the convergence  complexity of the pointwise iterates (see also   \cite[Theorem 6]{BMSZ21}) and the primal residual, where
the notation $\partial \theta (\m{x}) $ represents its sub-differential at $\m{x}$, and
$\C{N}_{\C{X}}(\m{x})$ denotes the normal cone of $\C{X}$ at $\m{x}$.
 \begin{theorem}  \label{b-10-1}
For any  $k>0$,  there exists an integer $t\leq k$  such that
 \begin{equation}\label{b-10-2}
   \big\|\mathbf{w}^{t+1} - \mathbf{w}^t \big\|^2_H  \le \frac{\varrho}{k+1}   \quad \textrm{ and } \quad
   \big\|\m{s}^t  \big\|^2_H \leq \frac{\varrho}{k+1}\frac{
   \big\| \operatorname{diag}(rA\tr A+Q, A\tr) \big\|^2_H}{ \gamma^2}  ,
\end{equation}
where $\m{s}^t\in \C{R}^n$ satisfies $A\tr \widetilde{\lambda}^t -\m{s}^t \in \partial \theta(\widetilde{\m{x}}^t)+
\C{N}_{\C{X}}(\widetilde{\m{x}}^t)$ and $\varrho=\frac{\gamma}{2-\gamma}\big\|\m{w}^*-\m{w}^0 \big\|^2_H$.
\end{theorem}
 Proof.
 Let  $k>0$ be a fixed constant and $t\leq k$ be a positive integer such that
\[{\big\|\m{w}^{t+1} - \m{w}^t\big\|^2_H} = \min \Big\{ {\big\|\m{w}^{l+1} - \m{w}^l\big\|^2_H}\mid l = 0, \ldots, k\Big\}.\]  Summing up (\ref{contract}) over $k=0,\cdots,\infty$ immediately gives
 \[
 \sum\limits_{k=0}^{\infty}\big\|\m{w}^k-\m{w}^{k+1} \big\|^2_H\leq\frac{\gamma }{2-\gamma}\big\|\m{w}^*-\m{w}^0 \big\|^2_H <\infty,
 \]
 which  further  shows
\begin{equation}\label{j-1}
{\big\|\m{w}^{t+1} - \m{w}^t  \big\|_H^2} \le \frac{\varrho}{k+1}.
\end{equation}

Now, it follows from (\ref{b-5})  that   by defining
\begin{equation}\label{right-eq2}
\m{s}^{t}=(rA\tr A+Q) (\widetilde{\m{x}}^t-\m{x}^t)+ A\tr (\widetilde{\lambda}^t-\lambda^t),
\end{equation}
we have
\[
\theta(\m{x}) - \theta(\widetilde{\m{x}}^t) {+}  (\widetilde{\m{x}}^t - \m{x}) \tr ( A   \tr \widetilde{\lambda}^t - \mathbf{s}^t)
\ge 0, \quad \forall \m{x}\in \C{X},
\]
which implies $A  \tr \widetilde{\lambda}^t - \mathbf{s} ^t \in \partial \theta(\widetilde{\m{x}}^t) +
\C{N}_{\C{X}}(\widetilde{\m{x}}^t)$.
By (\ref{right-eq2}) and (\ref{rela}) again, we have
\begin{eqnarray*}
{\big\| \mathbf{s}^t \big\|_H}
&=&  {\big\|\operatorname{diag}(rA\tr A+Q, A\tr) (\widetilde{\m{w}}^t-\m{w}^t)\big\|_H} =   \big\|\frac{1}{\gamma}\operatorname{diag}(rA\tr A+Q, A\tr) (\m{w}^{t+1}-\m{w}^t)\big\|_H \\
&\leq& \frac{1}{\gamma} \big\| \operatorname{diag}(rA\tr A+Q, A\tr) \big\|_H\big\|\m{w}^{t+1}-\m{w}^t\big\|_H,
\end{eqnarray*}
which together with (\ref{j-1}) ensures the right-hand-side inequality of  (\ref{b-10-2}).    $\hfill\blacksquare$

 The forthcoming remark suggests  that   our analysis of  Theorem \ref{b-10-1} is much easier than the analysis in e.g. \cite{HeY18}, although the   iteration-complexity results are consistent. The second remark shows that the pointwise iteration complexity can be still ensured  by regarding the relaxation factor $\gamma$ as some special  sequences.
\begin{remark}
Analogous to the technique \cite{HeY18} to establish the sublinear convergence rate in the pointwise sense, by adding the inequality (\ref{b-3}) with $\m{w}:=\widetilde{\m{w}}^{k+1}$, {i.e.,
 \[
 \theta(\widetilde{\m{x}}^{k+1})- \theta(\widetilde{\m{x}}^k) +\big\langle \widetilde{\m{w}}^{k+1}-\widetilde{\m{w}}^k, \C{J}(\widetilde{\m{w}}^{k+1})\big\rangle \geq   \big\langle\widetilde{\m{w}}^{k+1}-\widetilde{\m{w}}^k, H(\m{w}^k-\widetilde{\m{w}}^k) \big\rangle
 \]
 to the inequality (\ref{b-3}) at $(k+1)$-th iteration with $\m{w}:=\widetilde{\m{w}}^k$, i.e.,
  \[
 \theta(\widetilde{\m{x}}^k)- \theta(\widetilde{\m{x}}^{k+1}) +\big\langle \widetilde{\m{w}}^k-\widetilde{\m{w}}^{k+1}, \C{J}(\widetilde{\m{w}}^k)\big\rangle \geq   \big\langle\widetilde{\m{w}}^k-\widetilde{\m{w}}^{k+1}, H(\m{w}^{k+1}-\widetilde{\m{w}}^{k+1}) \big\rangle
  \]
  together} with the property (\ref{b-2}), we have
\[
  \big\langle \widetilde{\m{w}}^k -\widetilde{\m{w}}^{k+1}, H(\m{w}^k-\widetilde{\m{w}}^k +\widetilde{\m{w}}^{k+1}-\m{w}^{k+1}) \big\rangle  \geq 0.
\]
Then, adding the term $\big\|\m{w}^k-\widetilde{\m{w}}^k +\widetilde{\m{w}}^{k+1}-\m{w}^{k+1}  \big\|_H^2$ to both sides of the above inequality   and using $\m{w}\tr H\m{w}=\frac12\big\|\m{w}\big\|_{H+H\tr}^2$ and (\ref{rela}) will lead to
\begin{eqnarray}\label{cjb-1}
\frac12\big\|\m{w}^k-\widetilde{\m{w}}^k +\widetilde{\m{w}}^{k+1}-\m{w}^{k+1}  \big\|_{H+H\tr}^2
&\leq& \big\langle \m{w}^k -\m{w}^{k+1}, H(\m{w}^k-\widetilde{\m{w}}^k +\widetilde{\m{w}}^{k+1}-\m{w}^{k+1}) \big\rangle\nonumber\\
 &=&\gamma\big\langle \m{w}^k -\widetilde{\m{w}}^k, H(\m{w}^k-\widetilde{\m{w}}^k +\widetilde{\m{w}}^{k+1}-\m{w}^{k+1}) \big\rangle.\quad
\end{eqnarray}
So,  we   have from   (\ref{cjb-1}) and the identity ${\|a\|^2_H}- {\|b\|^2_H}= 2\langle a, H (a-b)\rangle - {\|a-b\|^2_H}$ that
\begin{eqnarray*}
&&\big\|\m{w}^k-\widetilde{\m{w}}^k\big\|_H^2 - \big\|\m{w}^{k+1}-\widetilde{\m{w}}^{k+1}\big\|_H^2\\
&=&2\big\langle \m{w}^k -\widetilde{\m{w}}^k, H(\m{w}^k-\widetilde{\m{w}}^k +\widetilde{\m{w}}^{k+1}-\m{w}^{k+1}) \big\rangle - \big\|\m{w}^k-\widetilde{\m{w}}^k +\widetilde{\m{w}}^{k+1}-\m{w}^{k+1}\big\|_H^2\\
&\geq&  \big\|\m{w}^k-\widetilde{\m{w}}^k +\widetilde{\m{w}}^{k+1}-\m{w}^{k+1}\big\|_{\widetilde{H}}^2,
\end{eqnarray*}
where
${\widetilde{H}}= \frac1\gamma(H+H\tr) -H=\frac{2-\gamma}{\gamma}H\succ \m{0}$.
By (\ref{rela}) again, the above inequality further shows
\[
\big\|\m{w}^k-\m{w}^{k+1}\big\|_H^2 \geq \big\| \m{w}^{k+1}-\m{w}^{k+2}\big\|_H^2,\quad \forall k\geq 0,
\]
which  using (\ref{contract})  gives
\begin{eqnarray*}
\frac{\gamma}{2-\gamma}\big\|\m{w}^{\kappa}-\m{w}^*\big\|^2_H\geq  \sum\limits_{k=\kappa}^{T+\kappa}
     \big\|\m{w}^k-\m{w}^{k+1}\big\|^2_H\geq (T+1)\big\|\m{w}^{T+\kappa}-\m{w}^{T+\kappa+1}\big\|^2_H.
\end{eqnarray*}
That is, $\big\|\m{w}^{T+\kappa}-\m{w}^{T+\kappa+1}\big\|^2_H\leq \frac{1}{T+1}\frac{\gamma}{2-\gamma}\big\|\m{w}^{\kappa}-\m{w}^*\big\|^2_H$ and it is consistent with the left-hand-side result of (\ref{b-10-2}). Especially, taking $\kappa=0$ is to obtain  $\big\|\m{w}^{T}-\m{w}^{T+1}\big\|^2_H\leq \frac{\varrho}{T+1}.$
\end{remark}

\begin{remark}\label{rem2}
Since the relaxation parameter $\gamma$ appears in the  iteration complexity results of  (\ref{b-10-2}), one may select a suitable fixed value  from the interval $(0,2)$ to do experiments or   update it adaptively by the following strategies    for any $k> 0$:
 \begin{enumerate}
  \item[(S1)]
  By updating   $\gamma^k=\frac{2k}{2k+c}$ where $c> 0$ is any real number,  we will have
  $\frac{\gamma^k}{2-\gamma^k}=\frac{k}{k+c}< 1$; or
  \item[(S2)]
   By updating   $\gamma^k=\frac{2k(k+2)}{2k^2+6k+3}$, we will obtain  $\frac{\gamma^k}{2-\gamma^k}=\frac{k(k+2)}{(k+1)(k+3)}<1;$
 \end{enumerate}
Then, it follows from (\ref{contract}) that
$
\big\|\m{w}^k-\m{w}^{k+1}\big\|^2_H \leq \big\|\m{w}^k-\m{w}^*\big\|^2_H - \big \|\m{w}^{k+1}-\m{w}^*\big\|^2_H
$
and then
 $
 \big\|\mathbf{w}^{t-1} - \mathbf{w}^t \big\|^2 \le \frac{1}{k} \big\|\m{w}^*-\m{w}^0 \big\|^2_H.
 $
 Meanwhile,   the convergence rates  of $\big\|\m{s}^t \big\|^2$ and   $\theta(\m{x}_t)- \theta(\m{x}^*) + \eta\big\|A\m{x}_t-b\big\|$ are  still $\C{O}({1}/{k}),$ {where $\eta$  is a constant.}

\begin{enumerate}
  \item[(S3)]
  Finally, at each step, we may choose a random number from 0 to 2. In this case, P-rALM still converges with a sublinear convergence rate.
 \end{enumerate}

\end{remark}

 \begin{theorem}  \label{b-10}
Let $T>0 $ and $\kappa\geq 0$. Then, for any $r>0 $ and $ Q\succ \m{0}$,  the sequence  $\{ \m{w}^k \}$   generated by P-rALM satisfies
 \begin{equation}\label{b-11}
  \theta(\m{x}_T)- \theta(\m{x}) + \big\langle  \m{w}_T-\m{w}, \C{J}(\m{w})\big\rangle\leq \frac{1}{2\gamma(T+1)}  \big\|\m{w}^\kappa-\m{w} \big\|^2_H,
  \ \forall \m{w} \in \C{M}.
\end{equation}
\end{theorem}
 Proof. The inequality (\ref{b-3}) and the relation (\ref{rela}) indicate
\begin{eqnarray}
 && \gamma \left[\theta(\m{x})- \theta(\widetilde{\m{x}}^k) +\big\langle \m{w}-\widetilde{\m{w}}^k, \C{J}(\m{w})\big\rangle\right] \geq \big\langle\widetilde{\m{w}}^k-\m{w}, H(\m{w}^{k+1}-\m{w}^k) \big\rangle\nonumber\\
 &= & \frac12\left\{ \big\|\widetilde{\m{w}}^k-\m{w}^k\big\|^2_H -\big\|\m{w}^{k+1}-\widetilde{\m{w}}^k\big\|^2_H + \big\|\m{w}^{k+1}- \m{w} \big\|^2_H -\big\|\m{w}^k- \m{w} \big\|^2_H \right\}, \label{key-11}
 \end{eqnarray}
in which the   equality uses the   identity
\begin{equation}\label{identabcd}
 \big\langle\m{p}-\m{q}, H (\m{u}-\m{v})\big\rangle =\frac12\Big\{ \big\|\m{p}-\m{v}\big\|^2_H - \big\|\m{p}-\m{u}\big\|^2_H +  \big\|\m{q}-\m{u}\big\|^2_H - \big\|\m{q}-\m{v}\big\|^2_H\Big\}
\end{equation}
with specifications
$
\m{p}:=\widetilde{\m{w}}^k,\m{q}=\m{w}, \m{u}=\m{w}^{k+1}$  and  $\m{v}:=\m{w}^k.
$
Note that
\begin{eqnarray*}
  \big\|\widetilde{\m{w}}^k-\m{w}^k\big\|^2_H -\big\|\m{w}^{k+1}-\widetilde{\m{w}}^k\big\|^2_H
 &= &\big\|\widetilde{\m{w}}^k-\m{w}^k\big\|^2_H -\big\|\m{w}^k-\widetilde{\m{w}}^k +\m{w}^{k+1}-\m{w}^k\big\|^2_H\\
 &= &\big\|\widetilde{\m{w}}^k-\m{w}^k\big\|^2_H -\big\|\m{w}^k-\widetilde{\m{w}}^k +\gamma(\widetilde{\m{w}}^k-\m{w}^k)\big\|^2_H\\
 &= & \gamma(2-\gamma)\big\|\widetilde{\m{w}}^k-\m{w}^k\big\|^2_H\geq 0.
 \end{eqnarray*}
Summing the   inequality (\ref{key-11}) over $k=\kappa,\kappa+1,\ldots,\kappa+T$, we have
\[
(T+1)\theta(\m{x})- \sum\limits_{k=\kappa}^{\kappa+T}\theta(\widetilde{\m{x}}^k) +\Big\langle (T+1)\m{w}-\sum\limits_{k=\kappa}^{\kappa+T}\widetilde{\m{w}}^k, \C{J}(\m{w})\Big\rangle  + \frac{1}{2\gamma}\big\|\m{w}-\m{w}^\kappa\big\|^2_H\geq 0,
\]
which, by the definition of $\m{w}_T$ and $\m{x}_T$, gives
\begin{equation}\label{b-12}
\frac{1}{T+1}\sum\limits_{k=\kappa}^{\kappa+T}\theta(\widetilde{\m{x}}^k) - \theta(\m{x})+ \big\langle \m{w}_T- \m{w}, \C{J}(\m{w})\big\rangle \leq \frac{1}{2\gamma(T+1)} \big\|\m{w}-\m{w}^\kappa\big\|^2_H.
\end{equation}
Because  $\theta$  is convex    and has the property
$
\theta(\m{x}_T)\leq\frac{1}{T+1}\sum_{k=\kappa}^{\kappa+T}\theta(\widetilde{\m{x}}^k),
$
   the inequality (\ref{b-11})  is obtained by plugging this property into (\ref{b-12}).$\hfill\blacksquare$

Theorem \ref{b-10} illustrates that the proposed P-rALM converges in a sublinear ergodic convergence rate. Furthermore, for any $\eta>0, $ by letting   $\Gamma_\eta=\{\lambda\mid \|\lambda\|\leq\eta\}$ and
\begin{equation}\label{b-12-1}
\gamma_\eta =\inf\limits_{\m{x}^*\in \C{X}}\sup\limits_{\lambda\in {\Gamma_\eta}} \big\| (\m{x}^*; \lambda)-(\m{x}^\kappa;\lambda^\kappa)\big\|_H^2,
\end{equation}
  we can  get the following tight result  whose proof is similar to that of \cite{bgc18,Shen20} and thus is omitted  here.  {We can see from Corollary \ref{b-15}    that  the residual  $\theta(\m{x}_T)- \theta(\m{x}^*) + \eta \|A\m{x}_T-b \|$ converges  with  the worst $\mathcal{O}(1/T)$ convergence rate.}
\begin{corollary}  \label{b-15}
 For any $\eta>0 $, let $\gamma_\eta$ be defined in (\ref{b-12-1}) and $\m{x}_T$ be defined in (\ref{erg-iterate}). Then,  the sequence  $\{ \m{w}^k \}$   generated by P-rALM satisfies
 \[
  \theta(\m{x}_T)- \theta(\m{x}^*) + \eta\big\|A\m{x}_T-b\big\|\leq \frac{\gamma_\eta}{2\gamma(T+1)},
  \ \forall \m{x}^*\in \C{X}.
\]
\end{corollary}

\subsection{Two special cases}\label{twocases}

In this section, we   discuss two interesting cases on  exploiting the Bregman distance and developing a linearized version of P-rALM, respectively.

 \textbf{Case 1}:  The double-penalty terms in P-rALM can be extended to the Bregman distance
\[
D_{\varphi}(\m{x},\m{x}^k)=\varphi(\m{x})-\varphi(\m{x}^k) -
\big\langle \nabla\varphi(\m{x}^k), \m{x}-\m{x}^k\big\rangle, \quad\forall \m{x}\in\C{X},
\]
where $\varphi$ is  a
 convex and    continuously differentiable function. For this case, by combining  the optimality condition (see e.g. \cite[lemma 2]{Zhalect13}) of the corresponding $\m{x}$-subproblem  and the previous {inequality} (\ref{b-6}), we can deduce
\begin{eqnarray*}
&& \big\langle\widetilde{\m{w}}^k-\m{w}^*, H(\m{w}^k-\widetilde{\m{w}}^k) \big\rangle - \big\langle\widetilde{\m{x}}^k-\m{x}^*, (r A\tr A +Q)(\m{x}^k-\widetilde{\m{x}}^k) \big\rangle
\\ &\geq & D_{\varphi}(\widetilde{\m{x}}^k,\m{x}^k) +
D_{\varphi}(\m{x}^*,\widetilde{\m{x}}^k) -
D_{\varphi}(\m{x}^*,\m{x}^k),
\end{eqnarray*}
that is
\begin{eqnarray}
&&\big\|\m{w}^*-\m{w}^k\big\|^2_H + 2 D_{\varphi}(\m{x}^*,\m{x}^k) -  \big\|\m{x}^*-\m{x}^k\big\|^2_{r A\tr A +Q}\nonumber\\
&\geq& \big\|\m{w}^*-\widetilde{\m{w}}^k\big\|^2_H + 2 D_{\varphi}(\m{x}^*,\widetilde{\m{x}}^k) -  \big\|\m{x}^*-\widetilde{\m{x}}^k\big\|^2_{r A\tr A +Q}\nonumber\\
&&+ \big\|\widetilde{\m{w}}^k- \m{w}^k\big\|^2_H + 2 D_{\varphi}(\widetilde{\m{x}}^k,\m{x}^k) - \big\|\widetilde{\m{x}}^k-\m{x}^k\big\|^2_{r A\tr A +Q}.\label{QA-ineq}
\end{eqnarray}
Denote $\tilde{H}=H-\operatorname{diag}(Q,\m{0})$. Under the assumption that {
\begin{equation}\label{Bre-bai}
D_{\varphi}(\m{x},\bar{\m{x}})\geq \frac{1}{2}
\big\{\big\| \m{x} -\bar{\m{x}}\big\|^2_{r A\tr A} -
\big\| \m{w} - \bar{\m{w}}\big\|^2_{\tilde{H}}  \big\}
\end{equation}
 for any $\m{w}, \bar{\m{w}}\in \C{M}$, the corresponding modified   P-rALM is convergent. When $\m{w}= \bar{\m{w}}$,   (\ref{Bre-bai}) holds obviously; when $\m{w}$ and $ \bar{\m{w}}$ are different points, we have $\big\| \m{w} - \bar{\m{w}}\big\|^2_{\tilde{H}}\geq 0$  and hence there exists a   constant $c>0$ such that $\big\| \m{w} - \bar{\m{w}}\big\|^2_{\tilde{H}}\leq \frac{c}{2}\big\| \m{x} -\bar{\m{x}}\big\|^2_{r A\tr A}.$ As a result, (\ref{Bre-bai}) reduces to the assumption on the kernel function $ \varphi$:
 \[
 \varphi(\m{x})-\varphi(\bar{\m{x}}) -
\big\langle \nabla\varphi(\bar{\m{x}}), \m{x}-\bar{\m{x}}\big\rangle\geq
\frac{1-c}{2}
\big\| \m{x} -\bar{\m{x}}\big\|^2_{r A\tr A}.
 \]
 If we take $\varphi=\frac12\big\|\m{x}\big\|^2_{r A\tr A}$, then $D_{\varphi}(\m{x},\bar{\m{x}})=  \frac12\big\|\m{x}-\bar{\m{x}}\big\|^2_{r A\tr A}$, showing that the assumption holds, and finally the   inequality (\ref{QA-ineq}) with simple algebra  will reduce to   (\ref{contract}). }

\textbf{Case 2}:  If the objective function $\theta(\m{x})$ is smooth and its gradient is  Lipschitz continuous with constant $L_\theta$, {which implies
 \begin{equation}\label{Lip-new}
 \theta(y) \leq  \theta(z) + \langle \nabla\theta(z), y-z\rangle+\frac{L_\theta}{2}\|y-z\|^2
 \end{equation}
for any $y,z\in\C{X}$,} then one may   update  the $ \widetilde{\m{x}}^k$-subproblem as the following
\begin{equation}\label{rem3-bj-1}
\widetilde{\m{x}}^k=\arg \min\limits_{\m{x} \in \C{X}} \langle\nabla\theta(\m{x}^k)-A\tr\lambda^{k},\m{x}\rangle   +\frac{r}{2}\big\|A(\m{x}-\m{x}^{k})\big\|^2 + \frac{1}{2}\big\|\m{x}-\m{x}^{k}\big\|^2_{Q},
\end{equation}
which, by taking $Q= \tau I - r A\tr A$ with $\tau> r \|A\tr A\|$, will become
$
\widetilde{\m{x}}^k=\C{P}_{\C{X}} \left(\m{x}^k-[\nabla\theta(\m{x}^k)-A\tr\lambda^k]/\tau \right).
$
  With the new update (\ref{rem3-bj-1}), a similar inequality to (\ref{contract}) can be also obtained. In fact, the first-order optimality conditions of (\ref{rem3-bj-1}) are $ \widetilde{\m{x}}^k\in\C{X}$ and
\[
  \big\langle\m{x}-\widetilde{\m{x}}^k, \nabla\theta(\m{x}^k) -A\tr \lambda^k+(rA\tr A+Q) (\widetilde{\m{x}}^k-\m{x}^k) \big\rangle\geq  0, ~\forall\m{x} \in\C{X},\\
\]
{which,  by using the convexity of $\theta$ and the inequality (\ref{Lip-new}) with $(y,z):=(\widetilde{\m{x}}^k, \m{x}^k)$, shows}
\begin{eqnarray*}
&& \big\langle\m{x}-\widetilde{\m{x}}^k,   A\tr \lambda^k+(rA\tr A+Q) (\m{x}^k-\widetilde{\m{x}}^k)\big\rangle  \\
 &\leq & \big\langle\m{x}-\widetilde{\m{x}}^k, \nabla\theta(\m{x}^k)  \big\rangle= \big\langle\m{x}-\m{x}^k+\m{x}^k-\widetilde{\m{x}}^k, \nabla\theta(\m{x}^k)  \big\rangle\\
 &\leq & \theta(\m{x})-\theta(\m{x}^k) +  \theta(\m{x}^k)- \theta(\widetilde{\m{x}}^k)+\frac{L_{\theta}}{2}\big\|\m{x}^k-
 \widetilde{\m{x}}^k\big\|^2.
\end{eqnarray*}
Rearrange it to get
\[
\theta(\m{x})-\theta(\widetilde{\m{x}}^k) + \langle\m{x}-\widetilde{\m{x}}^k,   -A\tr \lambda^k \rangle\geq  \langle\m{x}-\widetilde{\m{x}}^k,  (rA\tr A+Q) (\m{x}^k-\widetilde{\m{x}}^k) \rangle-\frac{L_{\theta}}{2}\|\m{x}^k-\widetilde{\m{x}}^k\|^2.
\]
Combine the last inequality with the   inequality (\ref{b-6}) to achieve
\[
\widetilde{\m{w}}^k\in \C{M},~  \theta(\m{x})- \theta(\widetilde{\m{x}}^k) +\left\langle \m{w}-\widetilde{\m{w}}^k, \C{J}(\m{w})\right\rangle \geq   \left\langle\m{w}-\widetilde{\m{w}}^k, H(\m{w}^k-\widetilde{\m{w}}^k) \right\rangle-\frac{L_{\theta}}{2}\|\m{x}^k-\widetilde{\m{x}}^k\|^2
\]
with $H$  given by (\ref{H}). Similar to the proof of Theorem \ref{Sec31-3}, we   deduce
\begin{equation}\label{rem3-bj-4}
   \big \|\m{w}^{k+1}-\m{w}^*\big\|^2_H\leq\big\|\m{w}^k-\m{w}^*\big\|^2_H -\frac{2-\gamma}{\gamma}
     \big\|\m{w}^k-\m{w}^{k+1}\big\|^2_{\tilde{H}}, \ \forall \m{w}^*\in \C{M}^*,
\end{equation}
where  $\tilde{H}=H-\operatorname{diag}\left(\frac{L_{\theta}}{2-\gamma}\m{I}, \m{0}\right)$. Clearly, if  $r>0$ and $Q\succ \frac{L_{\theta}}{2-\gamma}\m{I}$, then (\ref{rem3-bj-4}) implies that this linearized P-rALM   converges with the same convergence rate as P-rALM.

\section{\large Extensions of P-rALM for multi-block problem}\label{Sec3}
In this section, we {discuss two interesting extensions of P-rALM for solving a  multi-block separable convex minimization problem in the form of}
\begin{equation} \label{Sec4-Prob1}
\min  \Big\{\theta(\m{x}):=\sum\limits_{i=1}^{p}\theta_i(\m{x}_i)\Big|~ \sum\limits_{i=1}^{p}A_i\m{x}_i=b~ (\textrm{or} \geq b),~ \m{x}_i \in \C{X}_i \Big\},
\end{equation}
where      $ \theta_i: \C{R}^{n_i} \rightarrow
\C{R}, i=1,2,\cdots,p$, are    closed  proper convex functions; $\C{X}_i \subseteq\C{R}^{n_i}$ are   closed  convex   sets; $A_i\in\C{R}^{m\times n_i}$ and $b\in\C{R}^m$ are given data.  For this problem, we denote
\[
 \C{M}:=  \prod \limits_{i=1}^{p}\C{X}_i\times \Lambda,\ ~ \textrm{where}~
\Lambda:=\left\{\begin{array}{lll}
  \C{R}^m, & \textrm{ if } \sum_{i=1}^{p}A_i\m{x}_i=b, \\
 \C{R}^m_+, & \textrm{ if } \sum_{i=1}^{p}A_i\m{x}_i\geq b.
\end{array}\right.
\]
An extended primal-dual version of P-rALM   (denoted by \textbf{PD-rALM}) is described as follows.
\begin{flushleft}
\centering\fbox{
	\parbox{0.96\textwidth}{
\[
   \begin{array}{lll}
 \hspace*{-.1in}\textrm{\bf Initialize } (\m{x}_1^0,\ldots,\m{x}_p^0,\lambda^0), \textrm{ choose }
  \gamma\in(0,2), r_i>0,~ Q_i\succ\m{0} \textrm{ for } i=1,2,\ldots,p; \\
 \hspace*{-.1in} \textrm{\bf While } \textrm{the stopping criterion is not satisfied, }  \textrm{\bf do}\\
\hspace*{.1in}     \textrm{\bf For}\ i=1,2,\cdots,p, \textrm{\bf ~parallelly update}\\
\hspace*{.4in} \widetilde{\m{x}}_i^k=\arg \min\limits_{\m{x}_i \in \C{X}_i} \left\{\theta_i(\m{x}_i)- \langle \lambda^{k}, A_i\m{x}_i-b \rangle +\frac{r_i}{2}\big\|A_i(\m{x}_i-\m{x}_i^{k})\big\|^2 + \frac{1}{2}\big\|\m{x}_i-\m{x}_i^{k}\big\|^2_{Q_i}   \right\};\\
\hspace*{.1in} \textrm{\bf End for}\ \\
\hspace*{.1in}  \widetilde{\lambda}^k={\C{P}_{\Lambda}}\Big(\lambda^{k} -  \frac{1}{\sum_{j=1}^{p}\frac{1}{r_j}} \left[\sum_{i=1}^{p}A_i(2\widetilde{\m{x}}_i^k-\m{x}_i^k)-b\right]\Big);\\
\\
\hspace*{.1in}  \left(\begin{array}{c}
 \m{x}_1^{k+1}  \\ \vdots\\ \m{x}_p^{k+1}\\  \lambda^{k+1}
\end{array}\right) =  \left(\begin{array}{c}
 \m{x}_1^k  \\   \vdots\\ \m{x}_p^k \\ \lambda^k
\end{array}\right) + \gamma\left(\begin{array}{c}
\widetilde{\m{x}}_1^k- \m{x}_1^k \\ \vdots\\  \widetilde{\m{x}}_p^k- \m{x}_p^k\\  \widetilde{\lambda}^k- \lambda^k
\end{array}\right);\\
 \hspace*{-.1in} \textrm{\bf End while}
\end{array}
\]}}\end{flushleft}
Analogous to the analysis in Section \ref{Sec2}, the saddle-point  $\m{w}^*=(\m{x}_1^*;\cdots;\m{x}_p^*;\lambda^*)\in \C{M}$ of   (\ref{Sec4-Prob1}) will satisfy the previous
 $\textrm{VI}(\theta,\C{J},\C{M})$   but with new notations
\[
\m{w}=\left(\begin{array}{c}
 \m{x}  \\   \lambda
\end{array}\right)=
\left(\begin{array}{c}
 \m{x}_1  \\ \vdots\\  \m{x}_p\\   \lambda
\end{array}\right),\quad
\C{J}(\m{w})=\left(\begin{array}{c}
 -A_1\tr \lambda\\  \vdots\\ -A_p\tr \lambda\\ \sum_{i=1}^{p}A_i\m{x}_i-b
\end{array}\right).
\]

Compared to  the traditional ALM, the above PD-rALM is capable of exploiting the properties of each component objective function and could make the subproblems   much easier. Next,  we    briefly analyze the convergence of  this PD-rALM.
\begin{lemma} \label{1jx1-optimal}
The sequence  $\{ \m{w}^k \}$   generated by PD-rALM satisfies
\begin{equation}\label{jb-3}
\widetilde{\m{w}}^k\in \C{M},~ \theta(\m{x})- \theta(\widetilde{\m{x}}^k) +\big\langle \m{w}-\widetilde{\m{w}}^k, \C{J}(\m{w})\big\rangle \geq   \big\langle\m{w}-\widetilde{\m{w}}^k, H(\m{w}^k-\widetilde{\m{w}}^k) \big\rangle
\end{equation}
for any $\m{w}\in \C{M}$, where
\begin{equation}\label{j-H}
\widetilde{\m{w}}^k=\left(\begin{array}{c}
 \widetilde{\m{x}}^k  \\   \widetilde{\lambda}^k
\end{array}\right),~  H=\left[\begin{array}{cccccccc}
  r_1  A_1\tr A_1+Q_1 & \m{0}&  \cdots&  \m{0}& A_1\tr  \\
  \m{0} & r_2  A_2\tr A_2+Q_2&  \cdots&  \m{0}& A_2\tr  \\
  \vdots &  \vdots&\ddots&  \vdots& \vdots  \\
  \m{0} & \m{0}& \cdots&  r_p  A_p\tr A_p+Q_p& A_p\tr  \\
   A_1  & A_2&  \cdots&  A_p&  \sum\limits_{i=1}^{p}\frac{1}{r_i}  \mathbf{I}
\end{array}\right]
\end{equation}
is symmetric positive definite for any $r_i>0$ and $Q_i\succ \m{0}$. Moreover,   we have
\begin{eqnarray}\label{jb-4}
  \big \|\m{w}^{k+1}-\m{w}^*\big\|^2_H\leq\big\|\m{w}^k-\m{w}^*\big\|^2_H -\frac{2-\gamma}{\gamma}
     \big\|\m{w}^k-\m{w}^{k+1}\big\|^2_H, ~~ \forall \m{w}^*\in \C{M}^*.
\end{eqnarray}
\end{lemma}
 Proof.
First of all,    by the  first-order optimality condition of the $\m{x}_i$-subproblem ($i=1,2,\ldots,p$) in PD-rALM, we have $\widetilde{\m{x}}_i^k\in\C{X}_i$ and
\[
   \theta_i(\m{x}_i)- \theta_i(\widetilde{\m{x}}_i^k)  +  \left\langle\m{x}_i-\widetilde{\m{x}}_i^k, -A_i\tr \lambda^k+\big(r_iA_i\tr A_i+Q_i\big) (\widetilde{\m{x}}_i^k-\m{x}_i^k)\right\rangle\geq  0,   ~~ \forall \m{x}_i \in\C{X}_i,\]
 in other words,
\begin{eqnarray}\label{jb-5}
  \theta_i(\m{x}_i)- {\theta_i}(\widetilde{\m{x}}_i^k)  +  \left\langle\m{x}_i-\widetilde{\m{x}}_i^k, -A_i\tr \widetilde{\lambda}^k \right\rangle \geq    \left\langle\m{x}_i-\widetilde{\m{x}}_i^k,  \big(r_iA_i\tr A_i+Q_i\big) (\m{x}_i^k-\widetilde{\m{x}}_i^k)+ A_i\tr (\lambda^k-\widetilde{\lambda}^k)  \right\rangle.
 \end{eqnarray}
 Besides, it follows from the update of $\widetilde{\lambda}^k$ that $\widetilde{\lambda}^k\in \Lambda$ and
 \begin{equation}\label{jb-6}
  \Big\langle\lambda-\widetilde{\lambda}^k,  \sum\limits_{i=1}^{p}A_i\widetilde{\m{x}}_i^k-b  \Big\rangle{\geq}   \Big\langle\lambda-\widetilde{\lambda}^k,   \sum\limits_{i=1}^{p}A_i(\m{x}_i^k-\widetilde{\m{x}}_i^k) +\sum\limits_{j=1}^{p}\frac{1}{r_j}(\lambda^k-\widetilde{\lambda}^k)\Big\rangle,~\forall \lambda\in \Lambda.
\end{equation}
Finally,  combine the   inequalities (\ref{jb-5})-(\ref{jb-6}) together with the structure of $H$  and the {monotonicity} of $\C{J}(\m{w})$ to confirm the result in (\ref{jb-3}).

Note that the matrix
 $
H=\bar{H}+ \textrm{diag}(Q_1, \cdots,Q_p,\m{0})
$ where
 \begin{eqnarray*}
\bar{H}&=&\left[\begin{array}{ccccccc}
  r_1 A_1\tr A_1  &\cdots &\m{0}&  A_1\tr  \\
    \vdots  &\ddots&\vdots& \vdots  \\
     \m{0}  &\cdots & r_p A_p\tr A_p&  A_p\tr  \\
   A_1   &\cdots&  A_p&   \sum\limits_{i=1}^p\frac{1}{r_i} \mathbf{I}
\end{array}\right]\\
&=&\left[\begin{array}{ccccccc}
  r_1 A_1\tr A_1  &\cdots &\m{0}&  A_1\tr  \\
    \vdots  &\ddots&\vdots& \vdots  \\
     \m{0}  &\cdots & \m{0}&  \m{0} \\
   A_1   &\cdots&  \m{0}&   \frac{1}{r_1} \mathbf{I}
\end{array}\right]+\cdots+\left[\begin{array}{ccccccc}
  \m{0}  &\cdots &\m{0}&  \m{0} \\
    \vdots  &\ddots&\vdots& \vdots  \\
     \m{0}  &\cdots & r_p A_p\tr A_p&  A_p\tr  \\
   \m{0}   &\cdots&  A_p&    \frac{1}{r_p} \mathbf{I}
\end{array}\right]
\end{eqnarray*}
 \begin{eqnarray*}
&= & \left(\begin{array}{c}
  \sqrt{r_1} A_1\tr    \\
  \m{0}\\
  \vdots\\
  \m{0}\\
     \frac{1}{\sqrt{r_1}} \mathbf{I}
\end{array}\right)\left(\sqrt{r_1} A_1,\m{0},\ldots,\m{0},\frac{1}{\sqrt{r_1}} \mathbf{I}\right)+\ldots+ \left(\begin{array}{c}
 \m{0}   \\
   \vdots\\
  \m{0}\\
  \sqrt{r_p} A_p\tr\\
    \frac{1}{\sqrt{r_p}}\m{I}
\end{array}\right)\left(\m{0},\ldots,\m{0},\sqrt{r_p} A_p,\frac{1}{\sqrt{r_p}}\m{I}\right).
\end{eqnarray*}
For any $\m{w}=(\m{x};\lambda)\neq \m{0}$, we have
\[
\m{w}\tr H\m{w} = \sum\limits_{i=1}^p \Big\|\sqrt{r_i} A_i \m{x}_i+ \frac{1}{\sqrt{r_i}} \lambda\Big\|^2 +\sum\limits_{i=1}^p\big\| \m{x}_i \big\|_{Q_i}^2>0,
\]
and hence the matrix $H$ is symmetric positive definite. Similar to the proof of Theorem \ref{Sec31-3}, the inequality (\ref{jb-4}) can also be  obtained.
 $\hfill\blacksquare$

According to the previous  Lemma \ref{1jx1-optimal},  the global convergence of PD-rALM and its sublinear convergence rate    can be established  as   the rest parts of Section \ref{Sec2}.
{Motivated by the structure of   $H$ in (\ref{j-H}),   next   we present  a dual-primal update of PD-rALM, which can   also be regarded as a dual-primal extension of  P-rALM.

Let  $r_i>0, s_i>0$ and $ Q_i\succeq r_iA_i\tr A_i$ for   $i=1,2,\cdots,p$.  Consider the following block matrix
 \begin{equation}\label{H-new}
 H=\left[\begin{array}{cccccccc}
     Q_1 +s_1\mathbf{I}& \m{0}&  \cdots&  \m{0}& -A_1\tr  \\
  \m{0} & Q_2 +s_2\mathbf{I}&  \cdots&  \m{0}& -A_2\tr  \\
  \vdots &  \vdots&\ddots&  \vdots& \vdots  \\
  \m{0} & \m{0}& \cdots&  Q_p +s_p\mathbf{I}& -A_p\tr  \\
   -A_1  & -A_2&  \cdots&  -A_p&  \sum\limits_{i=1}^{p}\frac{1}{r_i}  \mathbf{I}
\end{array}\right].
\end{equation}
This new matrix $H$ is symmetric positive definite since
\begin{eqnarray*}
H\succeq{\footnotesize\underbrace{\left[\begin{array}{ccccccc}
  r_1 A_1\tr A_1  &\cdots &\m{0}&  -A_1\tr  \\
    \vdots  &\ddots&\vdots& \vdots  \\
     \m{0}  &\cdots & \m{0}&  \m{0} \\
   -A_1   &\cdots&  \m{0}&   \frac{1}{r_1} \mathbf{I}
\end{array}\right]+\cdots+\left[\begin{array}{ccccccc}
  \m{0}  &\cdots &\m{0}&  \m{0} \\
    \vdots  &\ddots&\vdots& \vdots  \\
     \m{0}  &\cdots & r_p A_p\tr A_p&  -A_p\tr  \\
   \m{0}   &\cdots&  -A_p&    \frac{1}{r_p} \mathbf{I}
\end{array}\right]}_{= {\footnotesize \left(\begin{array}{c}
  \sqrt{r_1} A_1\tr    \\
  \m{0}\\
  \vdots\\
  \m{0}\\
     -\frac{1}{\sqrt{r_1}} \mathbf{I}
\end{array}\right)\Big(\sqrt{r_1} A_1,\m{0},\ldots,\m{0},-\frac{1}{\sqrt{r_1}} \mathbf{I}\Big)+\ldots
+ \left(\begin{array}{c}
 \m{0}   \\
   \vdots\\
  \m{0}\\
  \sqrt{r_p} A_p\tr\\
    -\frac{1}{\sqrt{r_p}}\m{I}
\end{array}\right)\Big(\m{0},\ldots,\m{0},\sqrt{r_p} A_p,-\frac{1}{\sqrt{r_p}}\m{I}\Big)}}
+\left[\begin{array}{ccccccc}
   s_1\mathbf{I} &\cdots&  \m{0} &   \m{0}\\
   \vdots &\cdots&  \vdots &  \vdots \\
   \m{0} &\cdots&  s_p\mathbf{I} &  \m{0} \\
   \m{0}   &\ldots& \m{0} &  \m{0}
\end{array}\right],}
\end{eqnarray*}
and
\[
\m{w}\tr H\m{w} = \sum\limits_{i=1}^{p}\big\|\sqrt{r_i} A_i \m{x}_i- \frac{1}{\sqrt{r_i}} \lambda\big\|^2 +\sum\limits_{i=1}^{p}s_i\big\| \m{x}_i \big\|^2>0
\] for any $\m{w}=(\m{x};\lambda)\neq \m{0}$.
Substituting the above $H$ into (\ref{jb-3}), it is not difficult  to obtain the following dual-primal updates (denoted by {\bf DP-rALM}).
\begin{flushleft}
\centering\fbox{
	\parbox{0.97\textwidth}{
\[
  \begin{array}{lll}
 \hspace*{-.1in}\textrm{\bf Initialize } (\m{x}_1^0,\ldots,\m{x}_p^0,\lambda^0), \textrm{ choose }
  \gamma\in(0,2), r_i>0, s_i>0,~ Q_i\succ r_iA_i\tr A_i\textrm{ for } i=1,2,\ldots,p; \\
 \hspace*{-.1in} \textrm{\bf While } \textrm{the stopping criterion is not satisfied, }  \textrm{\bf do}\\
 \hspace*{.1in}  \widetilde{\lambda}^k=\C{P}_{\Lambda}\Big(\lambda^{k} - \frac{1}{\sum_{j=1}^{p}\frac{1}{r_j}} \Big[\sum_{i=1}^{p}A_i\m{x}_i^k-b\Big]\Big);\\
\hspace*{.1in}     \textrm{\bf For}\ i=1,2,\cdots,p, \textrm{\bf ~parallelly update}\\
\hspace*{.4in} \widetilde{\m{x}}_i^k=\arg \min\limits_{\m{x}_i \in \C{X}_i} \left\{\theta_i(\m{x}_i)- \langle  2\widetilde{\lambda}^k-\lambda^{k}, A_i\m{x}_i-b \rangle  + \frac{1}{2}\left\|\m{x}_i-\m{x}_i^{k}\right\|^2_{Q_i+ s_i\mathbf{I}}   \right\};\\
\hspace*{.1in} \textrm{\bf End for}\ \\
\\
\hspace*{.1in}      \left(\begin{array}{c}
 \m{x}_1^{k+1}  \\ \vdots\\ \m{x}_p^{k+1}\\  \lambda^{k+1}
\end{array}\right) =  \left(\begin{array}{c}
 \m{x}_1^k  \\   \vdots\\ \m{x}_p^k \\ \lambda^k
\end{array}\right) + \gamma\left(\begin{array}{c}
\widetilde{\m{x}}_1^k- \m{x}_1^k \\ \vdots\\  \widetilde{\m{x}}_p^k- \m{x}_p^k\\  \widetilde{\lambda}^k- \lambda^k
\end{array}\right);\\
 \hspace*{-.1in} \textrm{\bf End while}
\end{array}
\]}}\end{flushleft}

Although the condition of   $Q_i$ is stricter than that in   PD-rALM, similar results to Lemma \ref{1jx1-optimal} can be obtained.    {In fact, by the first-order optimality condition  of each  $\widetilde{\m{x}}_i^k$-subproblem ($i=1,2,\cdots,p$), we have} $\widetilde{\m{x}}_i^k\in \C{X}_i$ and
\begin{equation}\label{fin-1}
\theta_i(\m{x}_i) - \theta_i(\widetilde{\m{x}}_i^k)+\big\langle \m{x}_i-\widetilde{\m{x}}_i^k,-A\tr \widetilde{\lambda}^k \big\rangle\geq \big\langle \m{x}_i-\widetilde{\m{x}}_i^k,(Q_i+ s_i\mathbf{I})(\m{x}_i^k- \widetilde{\m{x}}_i^k )-A_i\tr(\lambda^{k}- \widetilde{\lambda}^k)\big\rangle
\end{equation}
for any $\m{x}_i \in\C{X}_i$. Meanwhile, the update of $\widetilde{\lambda}^k$ {implies $\widetilde{\lambda}^k\in \Lambda$ and}
\begin{equation}\label{fin-2}
\big\langle \lambda- \widetilde{\lambda}^k, \sum\limits_{i=1}^{p}A_i\widetilde{\m{x}}_i^k-b\big\rangle{\geq}\big\langle \lambda- \widetilde{\lambda}^k, - \sum\limits_{i=1}^{p}A_i(\m{x}_i^k-\widetilde{\m{x}}_i^k)+\sum\limits_{j=1}^{p}\frac{1}{r_j}
(\lambda^k-\widetilde{\lambda}^k)\big\rangle,~\forall \lambda\in \Lambda.
\end{equation}
Combining (\ref{fin-1}) and (\ref{fin-2}),  the previous inequalities (\ref{jb-3})
and (\ref{jb-4})   hold but with the matrix $H$ replaced by (\ref{H-new}). So,     DP-rALM   also converges  with a sublinear convergence rate.

}

\section{Numerical experiments}\label{Sec4}
In this section, we investigate the   performance of the proposed algorithms   for solving two popular optimization problems in machine learning. All the forthcoming experiments are   implemented in MATLAB R2019b
(64-bit) and performed on a PC
with Windows 10 operating system, with an Intel i7-8565U CPU and 16GB RAM.
\subsection{Linear support vector machine}
One fundamental function of machine learning is to make classification from a number of labeled training data. Suppose these   training data are   $ \{(x_i,y_i)\}_{i=1}^{m}$,  where $ x_i\in \C{R}^n$
are feature vectors   and $y_i\in\{-1,1\}$ are the labels  of   sample. If these   data formulate two disjoint convex hulls in $\C{R}^n$, then we can find a hyperplane $\{x\mid w\tr x +a=0\}$ to separate them because of the well-known strong separation theorem.
The linear support vector machine (abbreviated by SVM, see e.g. \cite{LiL15}) is to find the maximum margin hyperplane separating two classes of data as much as possible, which leads to the following optimization problem
\[
\min\limits_{w\in \C{R}^n, a\in \C{R}} \Big\{\frac{1}{2} \big\|w\big\|^2  \big|~ y_i\big(w\tr x_i +a\big)\geq 1, i=1,\cdots,m\Big\}.
\]
Introduce the following   notations
\[
u=\left(\begin{array}{c}
w \\ a
\end{array}\right),~ F=\left(\begin{array}{cc}
\m{I} & \m{0}\\ \m{0} & \m{0}
\end{array}\right),~ A =\left(\begin{array}{cc}
y_1(x\tr_1,1)  \\ \vdots \\ y_m(x\tr_m,1)
\end{array}\right)~~\textrm{and} ~~ b=\left(\begin{array}{c}
1\\ \vdots \\1
\end{array}\right)
\]
to reformulate the above SVM problem   as
\begin{equation}\label{Sec5-Pro}
\min\limits_{u\in \C{R}^{n+1}} \Big\{\frac{1}{2} \|Fu\|^2  \big|~ Au\geq b\Big\},
\end{equation}
which is clearly a special case of the problem (\ref{Sec1-Prob1}). {In fact, SVM has  been successfully applied in aero-engine fault diagnosis \cite{CYP15,HSYB05}. This subsection aims to test the numerical performance of our proposed methods for solving the problem (\ref{Sec5-Pro}).}
For example, applying the proposed P-rALM with $Q= \varrho \m{I} - r A\tr A$ to    (\ref{Sec5-Pro}), the resulting key iterations  are
\[\left\{\begin{split}
\widetilde{u}^k&=(F\tr F +\varrho \m{I})^{-1} (A\tr\lambda^k +\varrho u^k ),\\
   \widetilde{\lambda}^k&=\C{P}_{\C{R}^m_+}\left[\lambda^k - r (A(2\widetilde{u}^k-u^k)-b)\right].
\end{split}\right.
\]
An inverse operation is involved in the update of $\widetilde{u}^k$, but it is a fixed constant in each loop.

\begin{figure}[htbp]
 \begin{minipage}{1\textwidth}
 \def\figurename{\footnotesize Figure}
 \centering
 \resizebox{14.5cm}{7cm}
{
\includegraphics{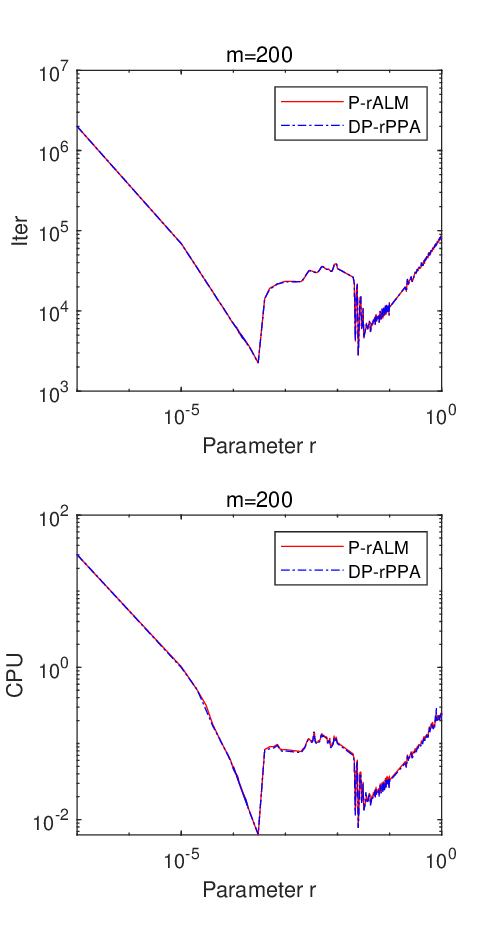}
\includegraphics{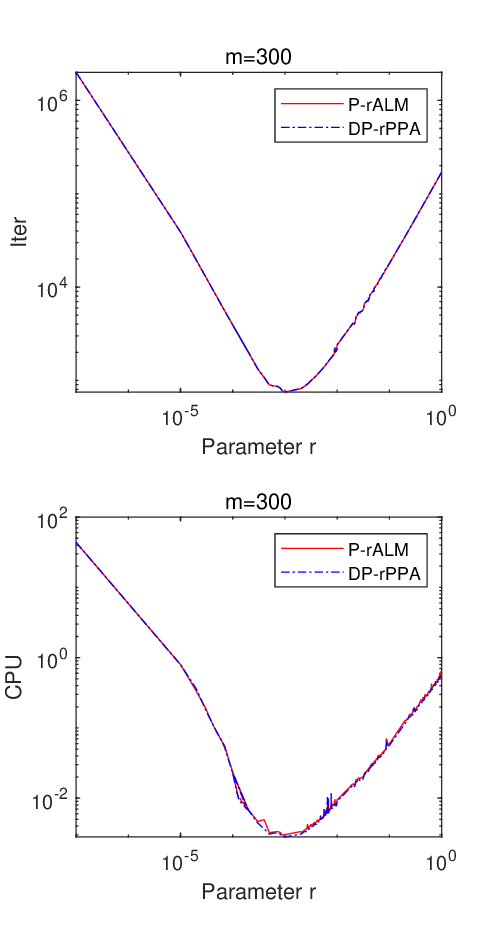}
\includegraphics{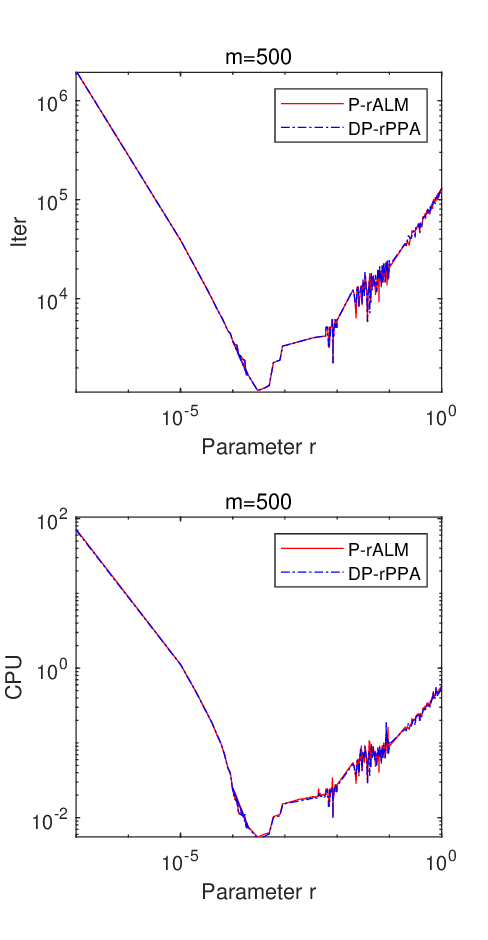}
}
\resizebox{14.5cm}{7cm}
{
\includegraphics{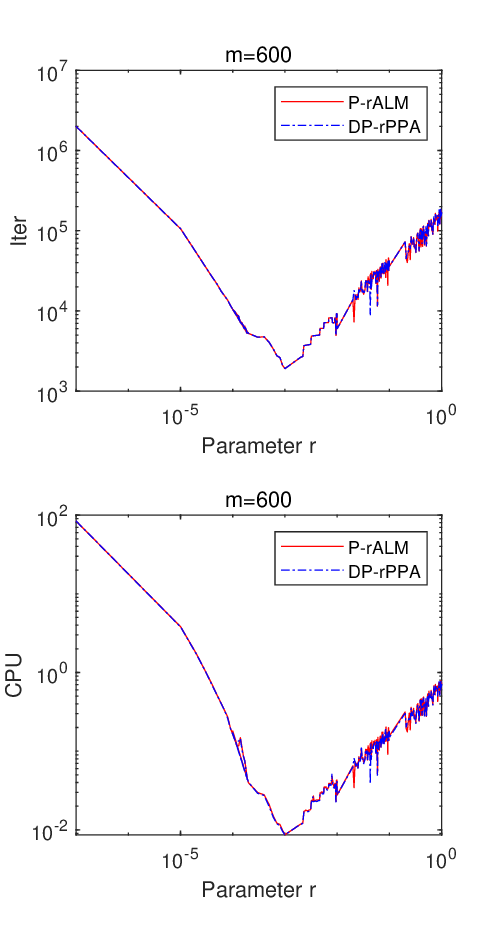}
\includegraphics{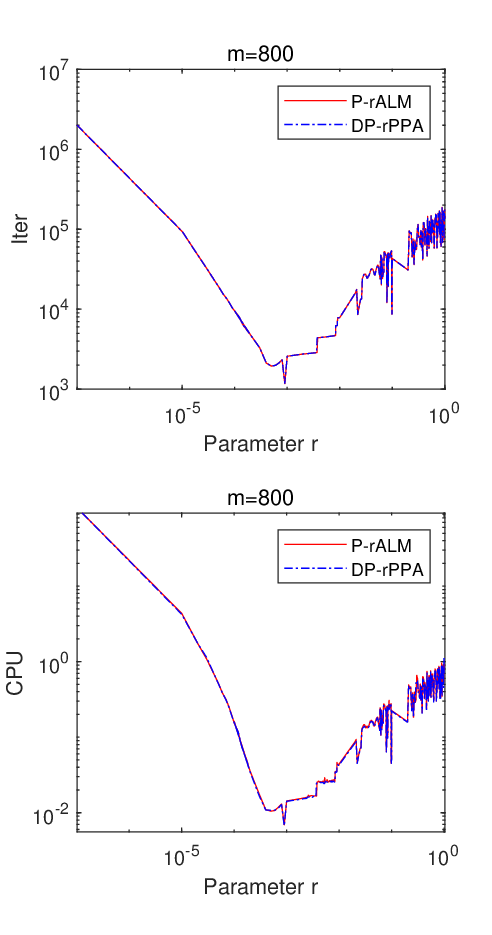}
\includegraphics{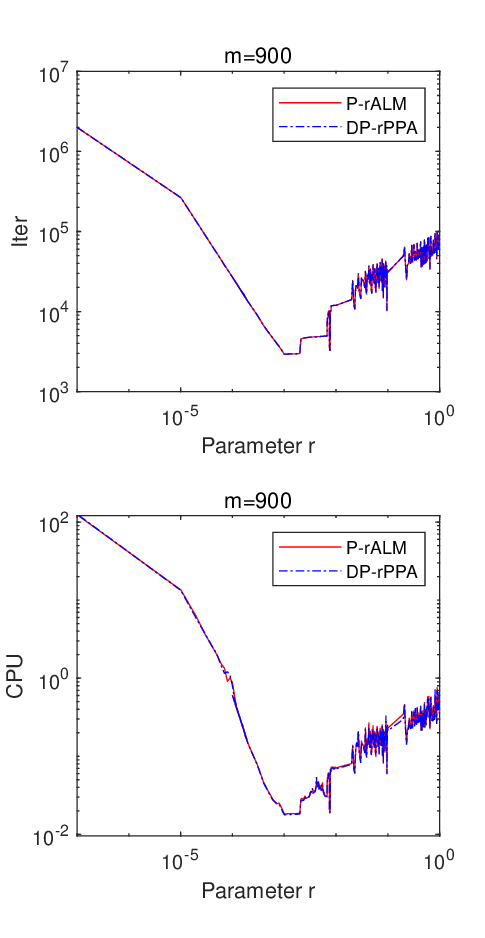}
}\vspace{-0.3cm}
\caption{\footnotesize Effect of  the  parameter $r$ on the performance of P-rALM and DP-rALM. }\label{testr}
   \end{minipage}
\end{figure}

The penalty parameter $r$ plays an important role in the performance of ALM.
To achieve a relatively better value of the penalty parameter $r$, we first test its effect   on the performance of our basic algorithms P-rALM and DP-rALM for solving the problem (\ref{Sec5-Pro}) with different   $m$. Throughout this subsection, the training data  are generated by random numbers
satisfying a normal distribution, and  the following stopping criterion
\begin{equation}\label{exper-1}
\textrm{Opt\_err}(k)= \max\left\{\|F\tr Fu^k -A\tr \lambda^k\|,\|\min(Au^k-b,\m{0})\|\right\} <tol
\end{equation}
is used to terminate P-rALM and DP-rALM under the maximal    iteration number  $2\times 10^{6}$. With  the same initial   points   $(u^0,\lambda^0)=(\verb"ones(3,1)",\verb"zeros(m,1)")$,  Figure  \ref{testr} shows the results about ``Iter" (the  number of iterations) and   ``CPU" (CPU time  in seconds) along with the increase of   $r$. {We fix  the tolerance   $tol=10^{-2}$  and the parameters    $(\varrho, \gamma)=(r (\rho(A\tr A)+0.1), 1.8)$ for  P-rALM, while   $(Q, s,\gamma)=(\varrho\m{I},10^{-3},1.8)$ for DP-rALM.}     We don't test  other values of $r$   since the reporting iteration numbers and CPU time are   worse than the results in Figure \ref{testr}. It can be seen from Figure \ref{testr} that both P-rALM and DP-rALM are competitive and sensitive to $r$. After checking  the reported  results {of both Iter and CPU} in Figure \ref{testr}, we find that $r=10^{-3}$  is   relatively reasonable    for   $m={600}$ to $900$, while $r=3\times10^{-4}$   for $m=200$ to {$500$}.

{To investigate the effect of the relaxation factor $\gamma\in(0,2)$ on the performance of    P-rALM and DP-rALM, Figure \ref{gamma-compa} presents  some comparative results of using different $\gamma$ under $tol=10^{-8}$, from which we can see that $\gamma=1.9$ performs relatively better  than others,  and it is set as the default   value in the following comparison  experiments. Besides, it can be seen from Figure \ref{gamma-compa} that  both P-rALM and DP-rALM enjoying a smaller  $\gamma\in(0,1)$   lead  to  much worse results.
\begin{figure}[htbp]
 \begin{minipage}{1\textwidth}
 \def\figurename{\footnotesize Figure}
 \centering
 \resizebox{14.5cm}{8cm}{
\includegraphics{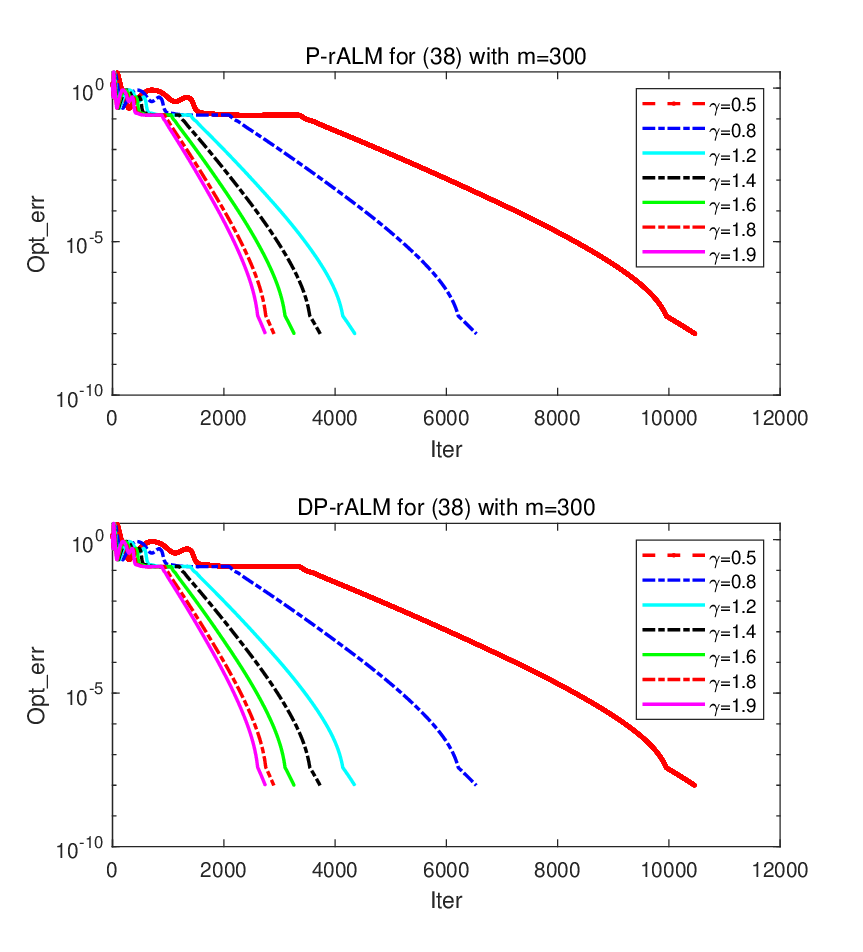}
\includegraphics{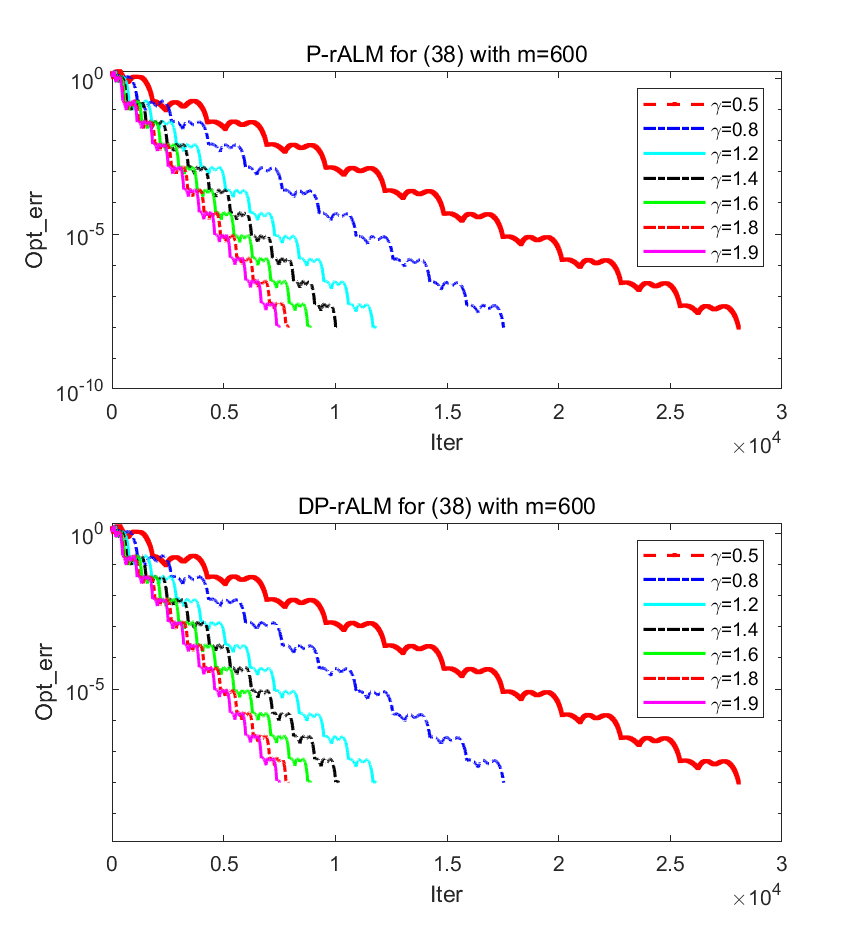}
\includegraphics{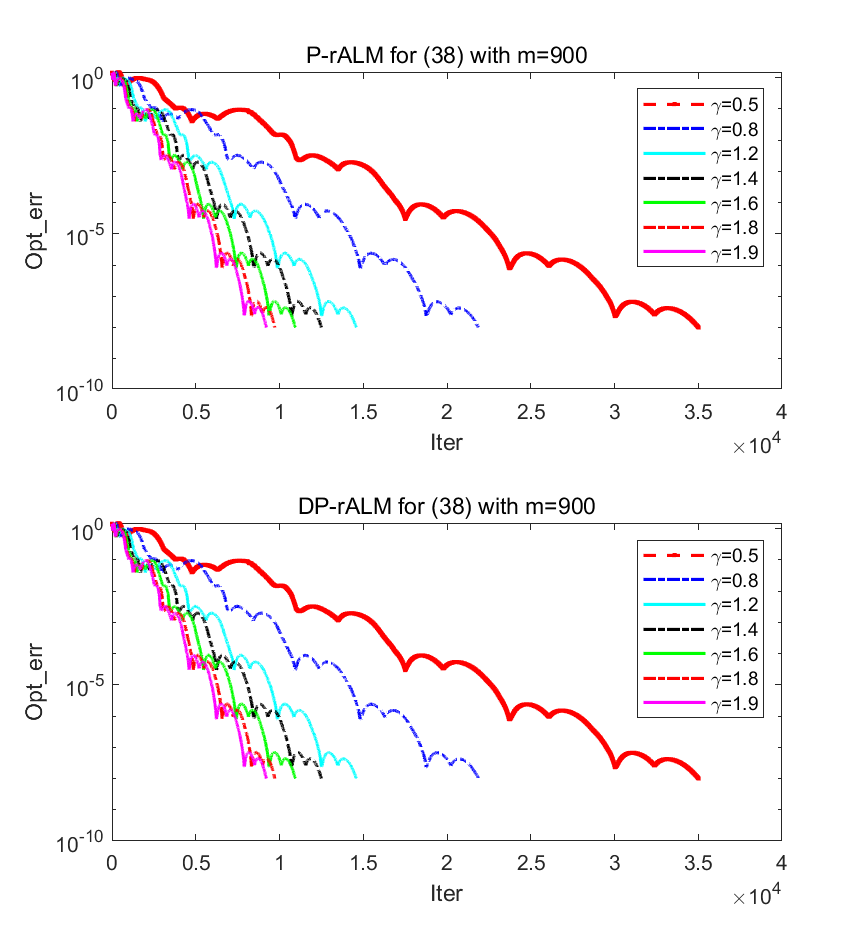}
}\vspace{-0.3cm}
\caption{\footnotesize Effect of  the relaxation parameter $\gamma$ on the performance of P-rALM and DP-rALM. }\label{gamma-compa}
   \end{minipage}
\end{figure}

A number of comparative experiments are further presented by   comparing   the following well-established   algorithms with our proposed algorithms P-rALM, DP-rALM, P-rALM with S1-S3,   DP-rALM with S1-S3 (the strategies S1-S3 are emphasized  in Remark \ref{rem2}), N-PDHG1 and N-PDHG2 using  the same relaxation step as   P-rALM:
 \begin{itemize}
   \item
    I-IDL-ALM \cite{HeXuYuj21} with parameters $(\tau,r)=(0.75,\beta \rho(A\tr A) +0.1)$    as the authors mentioned  but with     $\beta=0.004$  which
   performs better than the original setting $\beta=0.01$;
       \item
   C-PPA \cite{hy13} with $(\gamma, s) =(1.9, 1.01\rho(A\tr A)/r)$ where  we set  $r=m/8.5$;
   \item
   Generalized Primal-Dual Algorithm (G-PDA, \cite{HMan22}) with   primal-dual stepsize parameters \[
   \tau=c_1/\sqrt{(1-\alpha+\alpha^2)\rho(A\tr A)} \quad  \text{and}\quad \sigma=c_2/\sqrt{(1-\alpha+\alpha^2)\rho(A\tr A)}
   \]
    to satisfy the convergence condition $\frac{1}{\tau\sigma}>(1-\alpha+\alpha^2)\rho(A\tr A)$, and we set $(c_1,c_2,\alpha)=(2, 1/c_1-0.001,0.5)$ in the next experiments;
   \item
   G-PDHG \cite{Jiangzh23} with   the same setting for the parameters $(\tau, \sigma)$ as G-PDA  to satisfy the convergence condition $\frac{1}{\tau\sigma}>0.75\rho(A\tr A)$.
 \end{itemize}
 Here, we  emphasize that although the parameters of G-PDHG and  G-PDA are chosen in the same way,  their  frameworks are different. More specifically, the dual subproblem of G-PDA enjoys an expansion step with   inertial parameter $\alpha\in[0,1]$, while  G-PDA does not exploit  this step; G-PDA has a correction step for the dual-iterate, while G-PDHG has a correction step for the primal-iterate.
The  primal-dual methods  G-PDA, G-PDHG, N-PDHG1 and N-PDHG2 are used to solve the  saddle-point problem of (\ref{Sec5-Pro}):
$
\min_u\max_{\lambda\geq 0}  \frac{1}{2} \|Fu\|^2 -\langle \lambda, Au-b\rangle.
$
The}  parameters of our proposed algorithms use the tuned values as we mentioned before, and the parameter $c$  in  Remark \ref{rem2} is taken as 0.1.

\begin{table}
\centering
{\scriptsize{\begin{tabular}{cccccccccccccccc|} \hline
Size &\multicolumn{4}{c}{\textbf{P-rALM} }& &\multicolumn{4}{c}{\textbf{DP-rALM}} \\
\cline{2-5}\cline{7-10}m & Iter&  CPU  &\textrm{It\_err} &\textrm{Opt\_err} & & Iter  & CPU  &\textrm{It\_err}&\textrm{Opt\_err} \\
 \hline
300 &\emph{2747}   & \textbf{0.0163}    &  1.0250e-9 &    9.9916e-9    && \textbf{2746}	   & \emph{0.0169} &  1.0363e-9   &   9.9119e-9 \\
400 &\emph{11352}   &  0.1108 &  6.0521e-10 &    9.7526e-9    && \emph{11352}   &\textbf{0.0911} &  5.9308e-10   &   9.5964e-9 \\
500 &  \textbf{8115} &0.0942     & 2.3510e-10 &    9.9957e-9    &&  \emph{8116}  & \textbf{0.0671} & 2.3558e-10   &9.8216e-9 \\
600 & \emph{7518}  & \textbf{0.0567}   & 1.4976e-10 &    9.9693e-9    && \emph{7518}   &0.0650  & 1.4987e-10 &   9.9922e-9 \\
700 & \emph{13464}  & 0.1546    & 1.3101e-10 & 9.9949e-9    &&  \emph{13464}  & \textbf{0.1205} & 1.3062e-10   &   9.9744e-9 \\
800 &  15411 &  0.1324   &  6.0571e-11 &    9.8826e-9    && 15411   & 0.1268 &  6.0614e-11 &   9.8904e-9 \\
900 & \emph{9213}  &  0.0932   &  8.2811e-11 & 9.9745e-9    &&  \emph{9213}	  &  0.1169&  8.2823e-11   &   9.9762e-9 \\
\hline
Size &\multicolumn{4}{c}{\textbf{P-rALM with S1} }&&\multicolumn{4}{c}{\textbf{DP-rALM with S1}} \\
\cline{2-5}\cline{7-10}m & Iter&  CPU  &\textrm{It\_err} &\textrm{Opt\_err} & & Iter  & CPU  &\textrm{It\_err}&\textrm{Opt\_err} \\
  \hline
300 & 5229  & 0.0284    &  5.3706e-10 &    9.9946e-9    &&  5226  & 0.0302&  5.4623e-10   &   9.9787e-9 \\
400 & 21576&    0.1711 &  2.8475e-10 &    9.8869e-9    &&  21575  &0.1350	  &  2.7884e-10   &   9.9433e-9 \\
500 & 15427  & 0.1048    &  1.0254e-10 &    9.8898e-9    &&  15428	  & 0.1034 &  1.0238e-10	   &   9.9169e-9 \\
600 & 14028  &   0.0982  &  7.6785e-11 &    9.9603e-9    &&  14028  & 0.1158 & 7.6858e-11   &   9.9523e-9 \\
700 & 25591  & 0.1787    & 6.8756e-11 & 9.9926e-9    &&  25591  &0.1642	  &  6.8552e-11	  &   9.9717e-9 \\
800 &  29286 & 0.1856    &  3.2832e-11 &    9.9382e-9    && 29286   & 0.2233 &  3.2879e-11  &   9.9459e-9 \\
900 &  17508  & 0.1449    &  4.3011e-11 &    9.9947e-9    &&  17508  &0.1421  &  4.3017e-11   &   9.9966e-9 \\
\hline
Size &\multicolumn{4}{c}{\textbf{P-rALM with S2} }& &\multicolumn{4}{c}{\textbf{DP-rALM with S2}} \\
\cline{2-5}\cline{7-10}m & Iter&  CPU  &\textrm{It\_err} &\textrm{Opt\_err} & & Iter  & CPU  &\textrm{It\_err}&\textrm{Opt\_err} \\
  \hline
300 & 5237	  & 0.0329    &  5.3541e-10 &    9.9660e-9    &&  5234  &0.0284  &  5.4445e-10   &   9.9477e-9 \\
400 & 21585  & 0.1284    &  2.8481e-10 &    9.8627e-9    &&  21584  &0.1182  &  2.7889e-10   &   9.9193e-9 \\
500 & 15435  &   0.1122  &  1.0230e-10& 9.9962e-9    &&  15437  & 0.1154	 &  1.0255e-10  &   9.8377e-9 \\
600 & 14037  &     0.1453&  7.6807e-11	 &    9.8755e-9    && 14037   & 0.1030 &  7.6880e-11   &   9.8675e-9 \\
700 & 25600  &   0.2132  &  6.8760e-11 &    9.9934e-9    && 25600   & 0.1822 &  6.8555e-11 &   9.9726e-9 \\
800 &   29295 & 0.1888    &  3.2903e-11 &    9.9494e-9    &&   29295  &0.1926  & 3.2950e-11  &   9.9571e-9 \\
900 & 17517  &    0.1580	 &  4.2983e-11 &    9.9851e-9    &&17517    &0.1441  &  4.2990e-11   & 9.9870e-9 \\
\hline
Size &\multicolumn{4}{c}{\textbf{P-rALM with S3} }& &\multicolumn{4}{c}{\textbf{DP-rALM with S3}} \\
\cline{2-5}\cline{7-10}m & Iter&  CPU  &\textrm{It\_err} &\textrm{Opt\_err} & & Iter  & CPU  &\textrm{It\_err}&\textrm{Opt\_err} \\
 \hline
300 & 5245  & 0.0310  &  4.0011e-10 &   9.9876e-9    && 5257	   & 0.0448	 &  7.2148e-10   &  9.9437e-9 \\
400 &21370   &  0.1374   &  3.3819e-10 &    9.8724e-9    && 21593   &  0.1257&  4.7754e-10   &   9.8583e-9 \\
500 & 15386  & 0.1312    & 1.8824e-10 &    9.7337e-9    && 15483   & 0.1071 &  1.8188e-10 &   9.7578e-9 \\
600 & 14288  &    0.1476	 &  1.1117e-10	 &    9.9965e-9    &&  14160   & 0.1205 &  7.0168e-12  &   9.9993e-9 \\
700 &  25430  &  0.1951   &  4.7062e-11 &    9.9807e-9    &&  25717  & 0.1783	 &  1.3591e-10	 &   9.9731e-9 \\
800 &  29274 &   0.2409	  &  4.1498e-11  &    9.9694e-9    &&  29288  &  0.2682&  4.8997e-11 &   9.9513e-9 \\
900 & 17522	 &  0.2200   &  8.3521e-11 &  9.9658e-9    &&  17522  & 0.2013 &  6.3857e-11   &  9.9893e-9 \\
\hline
Size &\multicolumn{4}{c}{\textbf{I-IDL-ALM} }& &\multicolumn{4}{c}{\textbf{C-PPA}} \\
\cline{2-5}\cline{7-10}m & Iter&  CPU  &\textrm{It\_err} &\textrm{Opt\_err} & & Iter  & CPU  &\textrm{It\_err}&\textrm{Opt\_err} \\
 \hline
300 &  43099 &  0.2509  &  5.4987e-11 &    9.9979e-9    &&  4791	  &0.0347	  &  2.7731e-10   &   9.9823e-9 \\
400 & 166303  & 0.5714    &  4.0977e-11 & 9.9895e-9    &&  18504  &0.1404	  &  2.6511e-10   &   9.8792e-9 \\
500 &  88062 &  0.3580   &2.1486e-11 & 9.9728e-9    && 9840	   &0.0788  &  1.4402e-10  &   9.9608e-9 \\
600 &  35955 & 0.2108    & 4.2482e-11&    9.9628e-9    &&   17021  & 0.1465	 &  2.4382e-10 &   9.9742e-9 \\
700 &77014   &  0.4201   &  3.1677e-11 &    9.7316e-9    &&  \textbf{13047}  &0.1274  &  2.8562e-11	 &   9.9909e-9 \\
800 & 80489	  &   0.4837  &  2.0271e-11 &    9.9714e-9    &&  \emph{8872}  & \emph{0.0826} &  8.2374e-11 & 9.9742e-9 \\
900 & 47650  &  0.3991    &  1.9483e-11 & 9.9859e-9    &&   10383  &0.1382	  &  1.4948e-10   &   9.9613e-9 \\
\hline
Size &\multicolumn{4}{c}{\textbf{G-PDA} }& &\multicolumn{4}{c}{\textbf{G-PDHG}} \\
\cline{2-5}\cline{7-10}m & Iter&  CPU  &\textrm{It\_err} &\textrm{Opt\_err} & & Iter  & CPU  &\textrm{It\_err}&\textrm{Opt\_err} \\
 \hline
300 & 27323 &  0.1607   &  9.3168e-11 &   9.9999e-9    && 27590  &0.5680 &  9.3168e-11   &  9.9999e-9 \\
400 & 91678  &  0.3795   &  7.0041e-11 &    9.9984e-9    &&  87249  &1.0196  &  7.9698e-11   &   9.9924e-9 \\
500 & 47391  &  0.2808   & 1.5299e-11 &    9.9853e-9    &&  42990  &0.6013  &  2.3046e-11  &   9.9574e-9 \\
600 & 19717  &  0.1511   &  6.0752e-11 &    9.9740e-9    && 13448   & 0.2437 & 5.2588e-11  &   9.9894e-9 \\
700 & 29052	  & 0.2019    &  6.0593e-11 & 9.9939e-9    &&   31107  &  0.5585&  6.0509e-11	 &   9.9908e-9 \\
800 & 32069  &   0.2502  &  2.8833e-11 &    9.9958e-9    &&  31693  & 0.6043 & 5.6301e-11  &   9.9095e-9 \\
900 & 14741  &  0.1709   &  4.4148e-11 &    9.9775e-9    && 17512   & 0.3928 &  5.3862e-11   &   9.9896e-9 \\
\hline
Size &\multicolumn{4}{c}{\textbf{N-PDHG1} }& &\multicolumn{4}{c}{\textbf{N-PDHG2}} \\
\cline{2-5}\cline{7-10}m & Iter&  CPU  &\textrm{It\_err} &\textrm{Opt\_err} & & Iter  & CPU  &\textrm{It\_err}&\textrm{Opt\_err} \\
 \hline
300 & \emph{2747}  & 0.0231   &  1.0250e-9 &    9.9916e-9    &&  3410  &  0.0212	&  1.4921e-10   &   9.9875e-9 \\
400 & \emph{11352}  &   0.1095  &  6.0521e-10 &    9.7526e-9    &&  \textbf{11342}   & \emph{0.1076}&  5.9802e-10   &   9.7723e-9 \\
500 &   \textbf{8115} &  \emph{0.0691}   &  2.3510e-10 &    9.9957e-9    &&9590    & 0.0916	 &2.5845e-10   &   9.5993e-9 \\
600 &  \emph{7518}	  &    \emph{0.0570}	 & 1.4976e-10	 &    9.9693e-9    &&  \textbf{5720}  & 0.0595 &  9.2440e-11  &   9.9107e-9 \\
700 &  \emph{13464} &  \emph{0.1211}   & 1.3101e-10 &    9.9949e-9    &&13519    & 0.1520 &  1.2677e-10  &   9.9502e-9 \\
800 & 15411 &    0.1339 &  6.0571e-11	 &    9.8826e-9    &&  \textbf{6721}	  & \textbf{0.0538}	 & 3.1910e-11  &  9.9920e-9 \\
900 & \emph{9213 } &  \emph{0.0843}   &  8.2811e-11 &    9.9745e-9    &&  \textbf{7081}  & \textbf{0.0574} &  1.2773e-10   &   9.9970e-9 \\
\hline
\end{tabular}}
\caption{Comparative  results   of    different algorithms for  solving the problem  (\ref{Sec5-Pro}).}\label{TabProb1-1}}
\end{table}

Table \ref{TabProb1-1} {reports some} comparative   results    with different settings of total data number under $tol=10^{-8}$, where     $\textrm{Opt\_err}$ denotes  the final obtained   residual  defined in (\ref{exper-1}) and   $\textrm{It\_err}$ denotes the final obtained primal residual  defined as
$\textrm{It\_err}(k)= \|u^{k+1}-u^k\|.
$    The bold value in Table \ref{TabProb1-1}  denotes the smallest one for each problem and the italic value means    a  relatively smaller value. {We also set $n=2$ to visualize  the comparative  classification
results of $m\in\{300, 900\}$    as shown in
 Figures \ref{FigProb1-2}-\ref{FigProb1-3}.}  Figure \ref{Fig-cpma} depicts the convergence curves of  $\textrm{Opt\_err}(k)$   when applying these state-of-the-art algorithms to solve (\ref{Sec5-Pro}) with $m=900.$ {We have not depicted the curves about I-IDL-ALM in Figure \ref{Fig-cpma} since I-IDL-ALM performs much worse than others.} Figure  \ref{Fig-tol} further shows the comparison of different algorithms for solving (\ref{Sec5-Pro}) with $m=900$ under   lower and higher tolerance  errors.
 It can be seen from the  reported results   in Table \ref{TabProb1-1} and convergence curves   in Figures \ref{Fig-cpma}-\ref{Fig-tol} that:

\begin{figure}[htbp]
 \begin{minipage}{1\textwidth}
 \def\figurename{\footnotesize Figure}
 \centering
 \resizebox{14cm}{7.4cm}{\includegraphics{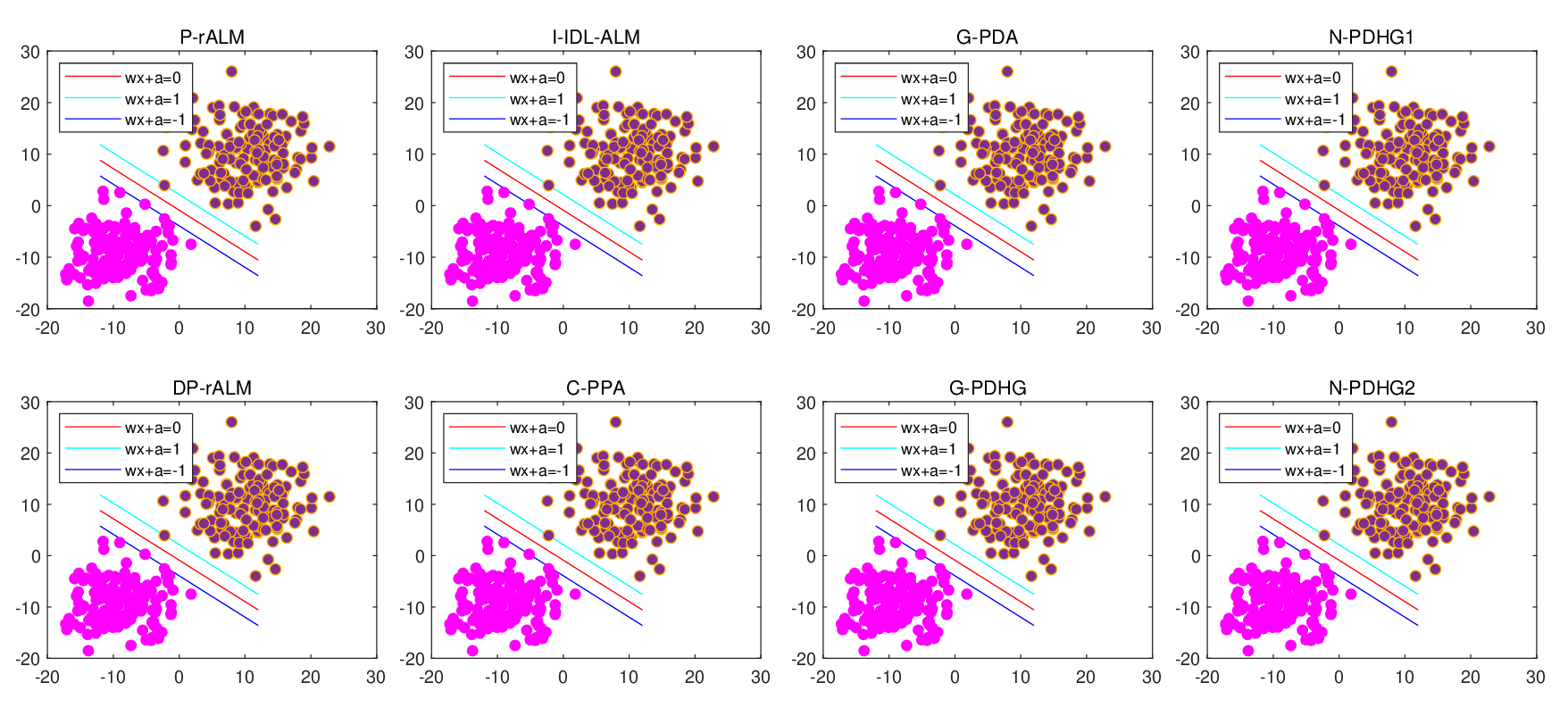}}\vspace{-0.3cm}
\caption{Classification results with  300 data  points by state-of-the-art algorithms.} \label{FigProb1-2}
 \centering
 \resizebox{14cm}{7.4cm}{\includegraphics{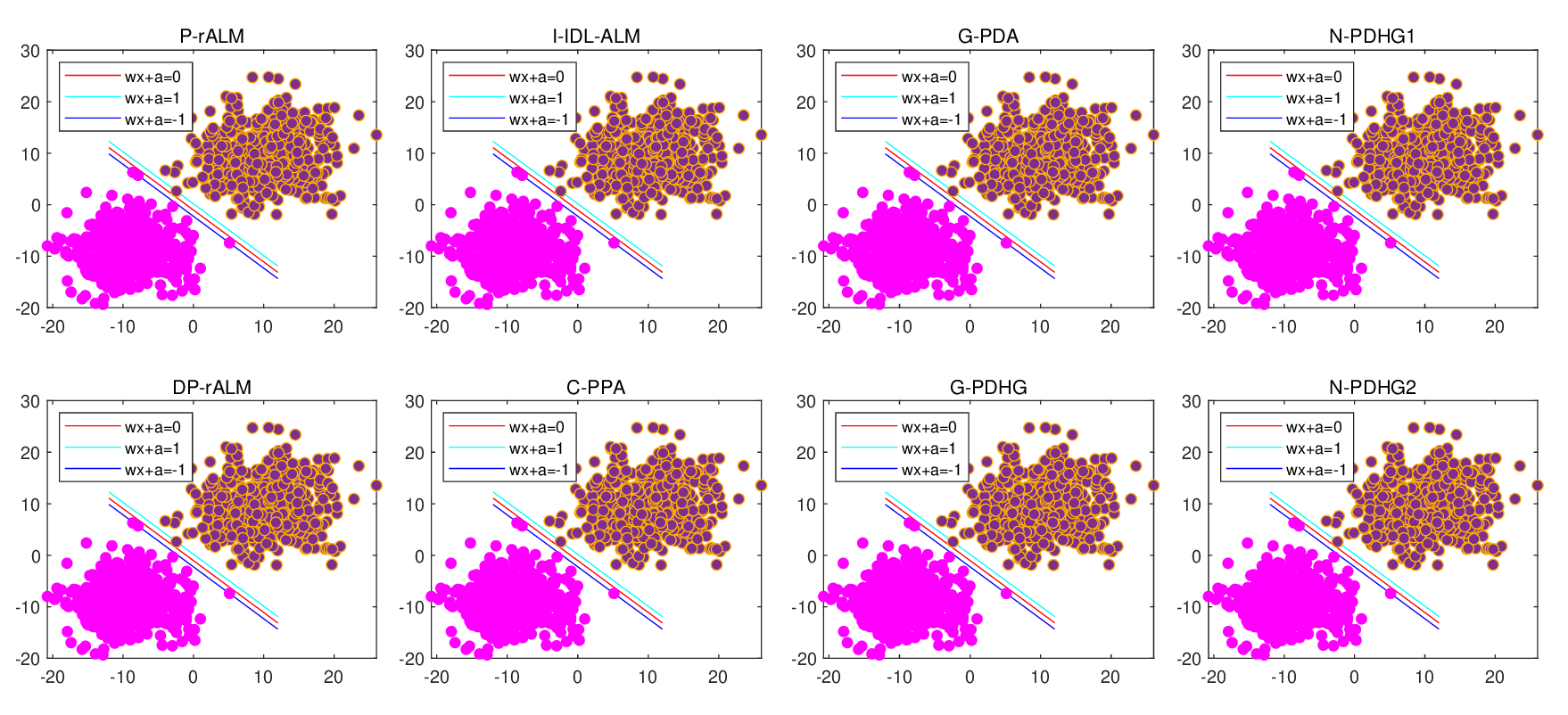}}\vspace{-0.3cm}
\caption{Classification results with 900 data points by state-of-the-art algorithms.} \label{FigProb1-3}
   \end{minipage}
\end{figure}

 \begin{itemize}
    \item
   {The  proposed methods  P-rALM and DP-rALM  are competitive and  perform significantly} better than each of them  with S1-S3, which suggests that using a large fixed relaxation factor is better than using an adaptive value.
   \item
   Both    P-rALM and DP-rALM     perform {significantly   better than  the recently  developed methods I-IDL-ALM, G-PDA and G-PDHG  for solving the problem (\ref{Sec5-Pro}) whether a} lower or higher tolerance is required.   {In addition, both N-PDHG1 and N-PDHG2 (see appendix)  are  competitive and perform also better than  I-IDL-ALM, G-PDA and G-PDHG, while N-PDHG2 seems more suitable for solving the large-scale SVM problem.}
   \item
  {Although C-PPA was proposed ten years ago,   it performs sometimes competitively with our P-rALM, which is perhaps due to the parameter choice  as  we suggested. However, in many cases,  this method performs   worse than our P-rALM and DP-rALM, which   can be observed  from  the reported    Iter and CPU.
     Compared to   all algorithms,}  I-IDL-ALM needs  relatively more iteration numbers to satisfy the stopping condition (\ref{exper-1}).
 \end{itemize}

 \begin{figure}[htbp]
 \begin{minipage}{1\textwidth}
 \def\figurename{\footnotesize Figure}
 \centering
 \resizebox{14cm}{5.2cm}{\includegraphics{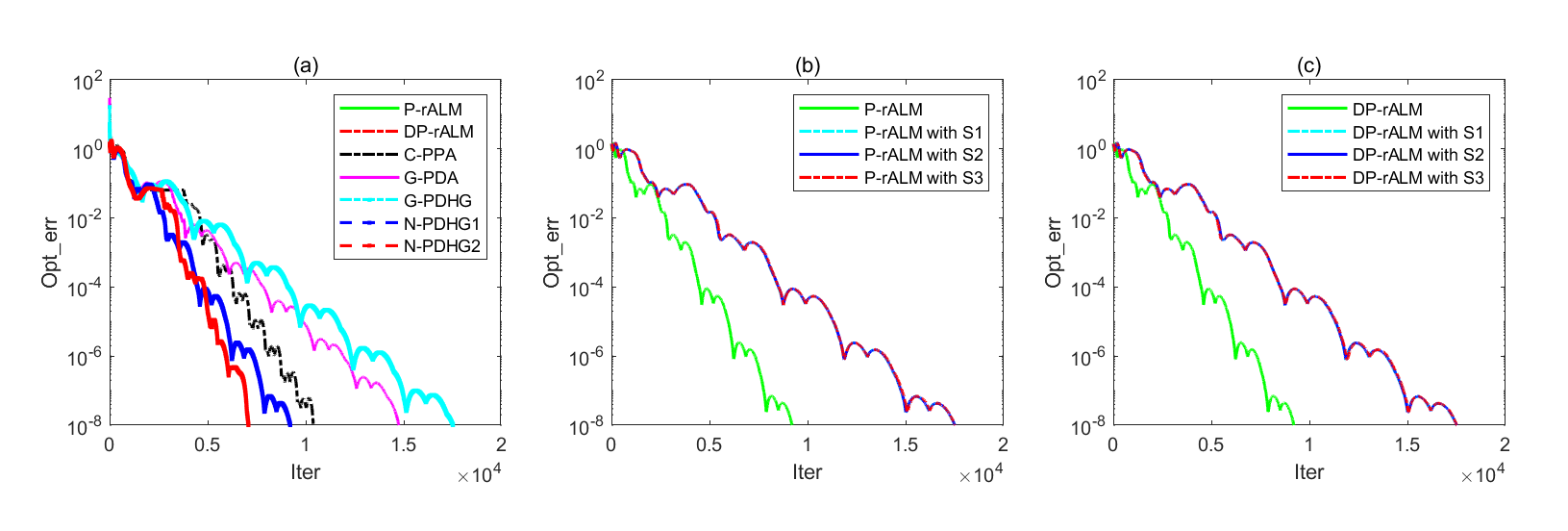}}\vspace{-0.4cm}
\caption{Comparison of Opt\_{err} versus Iter   by different algorithms for solving (\ref{Sec5-Pro}) with $m=900$.}
 \label{Fig-cpma}
 \centering
 \resizebox{14cm}{5.5cm}{\includegraphics{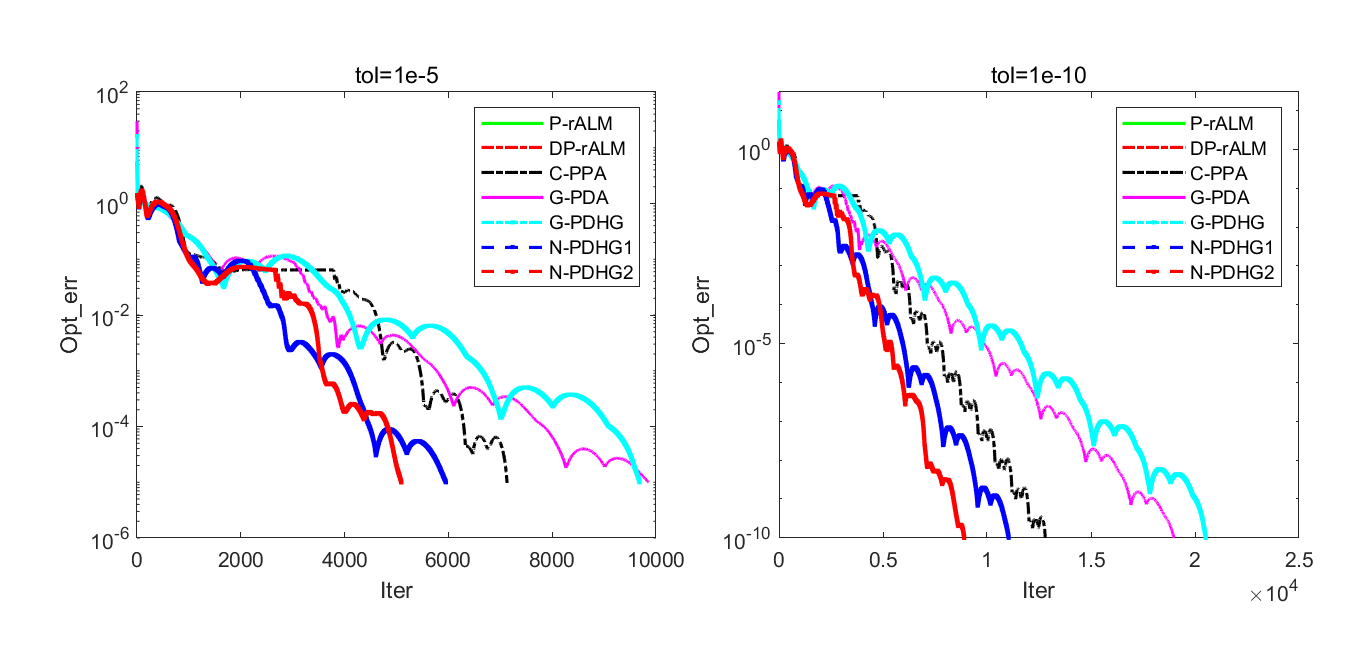}}\vspace{-0.4cm}
\caption{Comparison of different algorithms for solving (\ref{Sec5-Pro}) with $m=900$ under different $tol$.} \label{Fig-tol}
   \end{minipage}
\end{figure}

\subsection{Robust principal component analysis}
The Robust Principal Component Analysis (RPCA) was originally developed   by Candes et al. \cite{CLMW11}, which aims to decompose a data matrix $D\in \C{R}^{m\times n}$  into a   low-rank matrix $L$ and a sparse matrix $S$ containing outliers and corrupt data.  The principal components of $L$ are robust to the outliers and corrupt data in $S$. This
decomposition technique has   wide   applications in e.g.
video surveillance, face recognition,   latent semantic indexing, machine learning  and so forth, see e.g. \cite{ArgE07,SBNk19,CYP15,LLLBL17,Papa00}. Mathematically, the goal is to find $L$ and $S$   satisfying    the following separable nonconvex optimization problem
\[
\min\limits_{L,S\in \C{R}^{m\times n}} \big\{\textrm{rank}(L) + \|S\|_0 \mid  L+S=D\big\}.
\]
However, it is not a tractable optimization due to the non-convexity of  rank function  $\textrm{rank}(L)$ and the sparse norm $\|S\|_0.$  Similar to the technique to reformulate  the compressed sensing problem, most   researchers focus on   the following convex relaxation form:
\begin{equation}\label{Sec5-Pro2}
\min\limits_{L,S\in \C{R}^{m\times n}} \big\{ \|L\|_*  + \nu\|S\|_1 \mid  L+S=D\big\},
\end{equation}
where $\|\cdot\|_*$ denotes the nuclear norm of  a matrix (the sum of its singular values), $\|\cdot\|_1$
denotes the so-called $l_1$ norm of a matrix (the sum of its absolute values), and $\nu$ is a
positive weighting parameter that provides a trade-off between the sparse and low rank components and  usually $\nu=1/\sqrt{\max(m,n)}$.  Clearly, the problem (\ref{Sec5-Pro2}) is a special case of the previous   model (\ref{Sec4-Prob1}) and   hence {the     PD-rALM and DP-rALM presented  in section \ref{Sec3}} can be applied to solve it. For example,   applying PD-rALM with $Q_i=\varrho_i\m{I}$ for $i=1,2$, the resulting key iterations are
\[  \left\{\begin{split}
\widetilde{L}^k&=\arg\min\limits_{L\in \C{R}^{m\times n}} \|L\|_* + \frac{r_1+\varrho_1}{2}\Big\|L- L^k -\frac{\Lambda^k}{r_1+\varrho_1} \Big\|_F^2,\\
   \widetilde{S}^k&=\arg\min\limits_{S\in \C{R}^{m\times n}} \nu\|S\|_1 + \frac{r_2+\varrho_2}{2}\Big\|S- S^k -\frac{\Lambda^k}{r_2+\varrho_2} \Big\|_F^2\\ &= \textrm{Shrink}\Big(S^k +\frac{\Lambda^k}{r_2+\varrho_2}, \frac{\nu}{r_2+\varrho_2} \Big),\\
\widetilde{\Lambda}^k&=\Lambda^k  - \frac{1}{1/r_1 +1/r_2}  \Big[(2\widetilde{L}^k-L^k) + (2\widetilde{S}^k-S^k) -D\Big],
\end{split}\right.
\]
where $\textrm{Shrink}(\cdot,\cdot)$ denotes the soft shrinkage operator (see e.g. \cite{TY11}). And the $L$-subproblem admits the following explicit solution
\[
\widetilde{L}^k=U^k\operatorname{diag}\Big(\max\Big\{\sigma_i^k-\frac{1}{r_1
+\varrho_1},0\Big\}\Big)
(V^k)^T,
\]
where  $U^k\in \C{R}^{m\times r}$ and $V^k\in \C{R}^{n\times r}$ are obtained by the    singular value decomposition:
$L^k +\frac{\Lambda^k}{r_1+\varrho_1}=U^k\Sigma^k(V^k)^T
$
with
$
{\Sigma}^k=\operatorname{diag}\left(\sigma_1^k,\sigma_1^k,\cdots,\sigma_r^k\right).
$


{   In what follows, by comparing  with several well-established algorithms  we   test the performance of the preliminary algorithms   PD-rALM and DP-rALM} for solving the problem (\ref{Sec5-Pro2}) with  Yale B database\footnote{The database can be downloaded at http://vision.ucsd.edu/$\sim$iskwak/ExtYaleDatabase/ExtYaleB.html.} which consists of cropped and aligned images   of 38 individuals   under 9 poses and 64 lighting conditions.  The penalty parameter of the standard ADMM is  fixed as  $\frac{mn}{4\|D\|_1}$ according to \cite[Page 109]{SBNk19}. {Similar to  the parameter  choice  in ADMM, we set $r_1=r_2, Q_1=Q_2=\varrho\m{I}$ with $(r_1,\gamma, \varrho)=(\frac{mn}{5\|D\|_1},1.75,10^{-6})$ for  PD-rALM, while we set $s_1=s_2$  with $(\varrho,s_1)=(r_1(1+10^{-3}), 10^{-4})$ for DP-rALM.
  The existing   DP-BALM \cite{xushengjie23}, PDHG \cite{chpo11}, G-PDHG and G-PDA are used to solve the corresponding saddle-point problem
$
\min\limits_{L, S}\max\limits_{Z}
\{\|L\|_* +\nu \|S\|_1 + \langle L+S-D, Z\rangle\}.
$ The  parameters of DP-BALM  are fixed as $(\alpha, \beta, \delta)=( 1.9,100, 10^{-3})$, both G-PDA and  G-PDHG use  the suggested values as mentioned in the second experiments \cite{xushengjie23} since the involved  parameters  $(\tau,\sigma)$ are restricted by the same condition $\tau\sigma \|KK\tr\|<4/3$ and here $K=[\m{I};\m{I}]$. Due to the recent work \cite{BUGi23}, we set the inertial  parameter $\theta=0.8$ for PDHG. }
All mentioned  algorithms {are terminated when
  the following criteria
 \[
     \textrm{RelChg(k)}:=\frac{\big\|L^{k+1}-L^k\big\|_F+\big\|S^{k+1}-S^k\big\|_F}{
                       \big\|L^k\big\|_F+\big\|S^k\big\|_F+1}<\epsilon_1,\ ~
   \textrm{Res(k)}:=\frac{\|D - L^{k+1}-S^{k+1}\|_F}{\|D\|_F}<\epsilon_2
     \]
are satisfied with}   the same initial feasible points  $(\Lambda^0,S^0)=(L^0,D-L^0)$, where $(\epsilon_1,\epsilon_2)$ are given tolerances and  $L^0$ is obtained by the truncated singular value decomposition:
\[
  L^0=\verb"F(:,1:l)Sigma(1:l,1:l)N(:,1:l)" ~\textrm{where}~ \verb"[F,Sigma,N]=svd(D,'econ');"  \verb"l=3."
\]

\begin{table}[h]
{\centering
\begin{tabular}{ll l l   l l l}
\hline
$(\epsilon_1,\epsilon_2)$& Methods &  Iter & $\operatorname{rank}(L)$ & CPU  & RelChg & Res  \\
\hline
& PD-rALM  &\emph{194}   &37 & \textbf{17.6648} & 9.9673e-5 & 6.8302e-6\\
& DP-rALM  &\textbf{192} &37 & \emph{18.4554}   &9.9661e-5  & 6.8722e-6 \\
$(10^{-4}, 10^{-5})$& ADMM  &254&  31& 21.0341 & 9.8391e-5  &4.0032e-6 \\
& DP-BALM   &225   & 31 & 21.7985& 5.9275e-5  &9.9039e-6 \\
&    PDHG   & 280  & 31 & 23.0072 & 1.9046e-5  &9.9674e-6	\\
&   G-PDHG  & 271 &  31 & 22.5811 & 2.0347e-5   &9.9877e-6 \\
&   G-PDA   &271  &  31 & 22.0591 & 1.9011e-5  &9.9877e-6 \\
\hline
&   PD-rALM & \textbf{343}&  31  & \textbf{30.8044}& 9.8724e-6 & 7.4143e-7 \\
&  DP-rALM  &\emph{360}	    & 31   &\emph{32.4381} & 9.9899e-6& 5.6868e-7 \\
$(10^{-5}, 10^{-6})$&   ADMM  & 395	   &31    &38.1146&9.8640e-6& 4.7239e-7\\
&   DP-BALM   & 469   &  31  & 46.5738 &8.7892e-6 & 9.9287e-7 \\
&    PDHG   &  516  & 31   & 53.5624&   1.8649e-6 & 9.9971e-7 \\
&   G-PDHG  &  508  & 31   & 49.6827&2.0340e-6    & 9.9375e-7 \\
&   G-PDA   &  508  & 31   &48.8102& 1.9017e-6   & 9.9375e-7 \\
\hline
&   PD-rALM & \textbf{582}	&  31&\textbf{53.7565}          &9.9349e-7    &9.5168e-8  \\
&  DP-rALM  & \emph{605}&  31 & 56.6440  & 9.9954e-7   &8.6355e-8  \\
$(10^{-6}, 10^{-7})$&   ADMM  & 619 & 31&57.4334 & 9.9702e-7  & 4.4109e-8	 \\
&   DP-BALM   & 864   & 31  &85.6948 & 9.9842e-7  & 7.0061e-8 \\
&    PDHG   & 988  & 31 &93.3656 &1.2053e-7 & 9.9848e-8 \\
&   G-PDHG  & 976   & 31 &89.9472  &  1.4532e-7 & 9.9761e-8 \\
&   G-PDA   & 976   & 31 & 86.8632 &  1.2596e-7  & 9.9761e-8 \\
\hline
\end{tabular} \vspace{-0.2cm}
\caption{Comparative results of  state-of-the-art algorithms under different tolerances.} }\label{TabProb-b}\vspace{-0.3cm}
\end{table}

\begin{figure}[http]
 \begin{minipage}{1\textwidth}
 \def\figurename{\footnotesize Figure}
 \centering
 \resizebox{14.8cm}{5cm}{\includegraphics{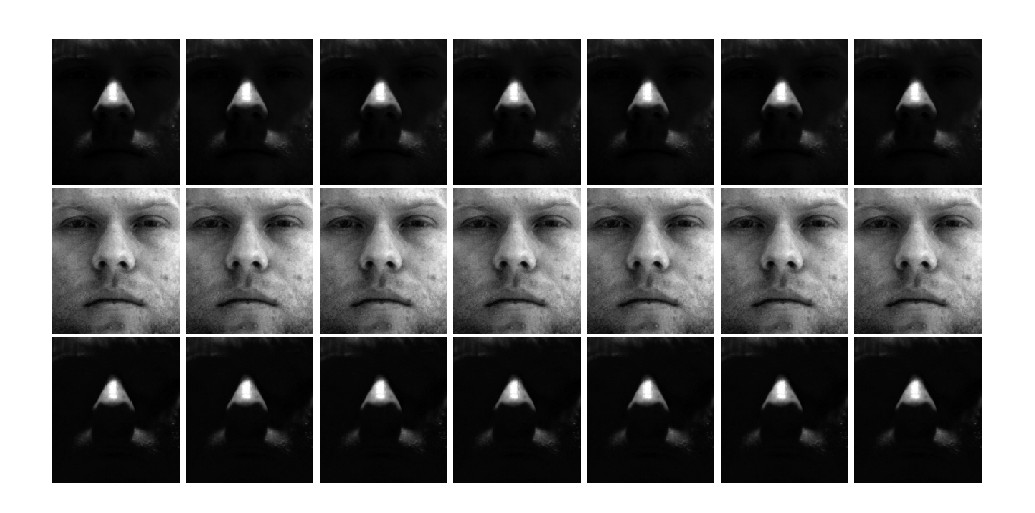}}\vspace{-0.5cm}
 \resizebox{14.8cm}{5cm}{\includegraphics{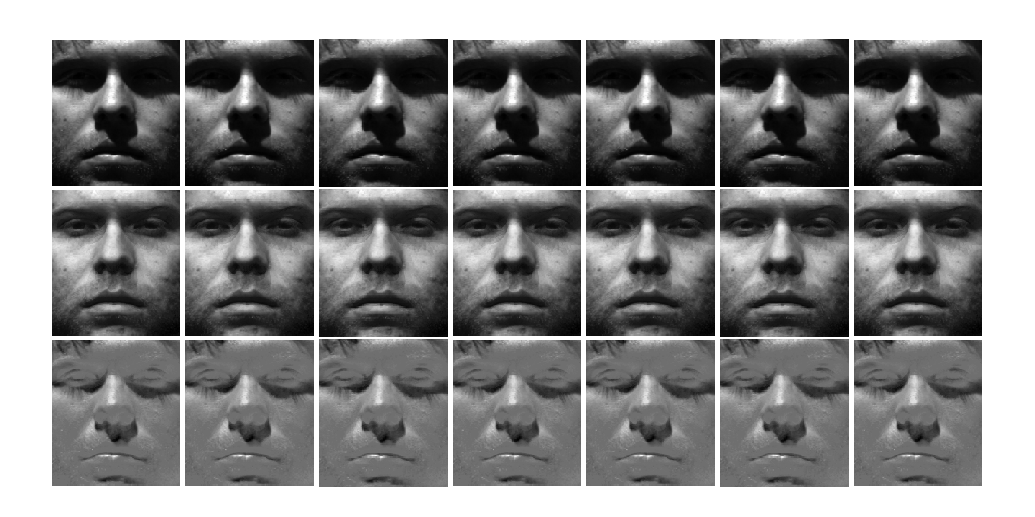}}\vspace{-0.5cm}
 \resizebox{14.8cm}{5cm}{\includegraphics{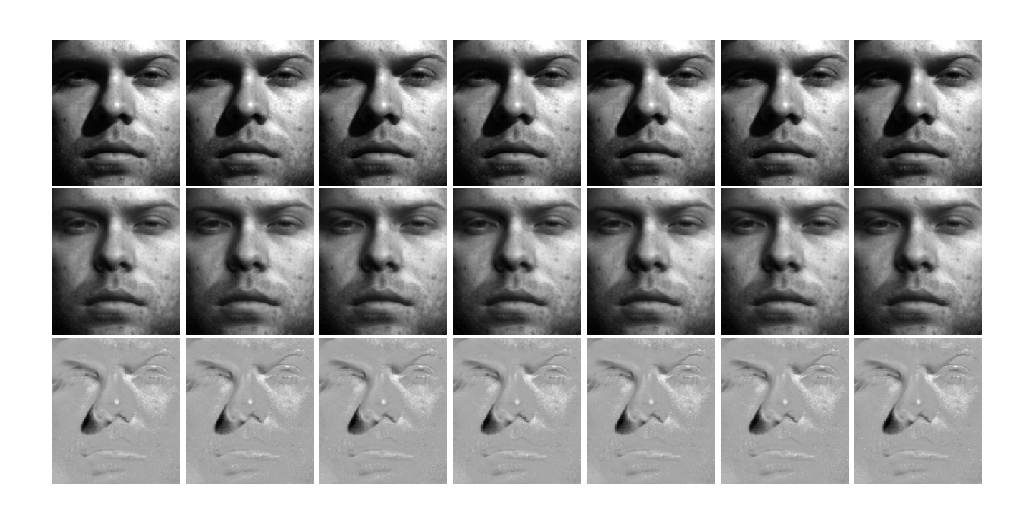}} \vspace{-0.3cm}
\caption{\footnotesize Output of different algorithms for the 4th(rows 1-3), 18th(rows 4-6) and 46th(rows 7-9) images  in the Yale B database. From left to right: PD-rALM, DP-rALM, ADMM, DP-BALM, PDHG, G-PDHG,   G-PDA, respectively.}\label{Fig4-2}
   \end{minipage}
\end{figure}

\begin{figure}[htbp]
 \begin{minipage}{1\textwidth}
 \def\figurename{\footnotesize Figure}
 \centering
 \resizebox{15.8cm}{5.2cm}{\includegraphics{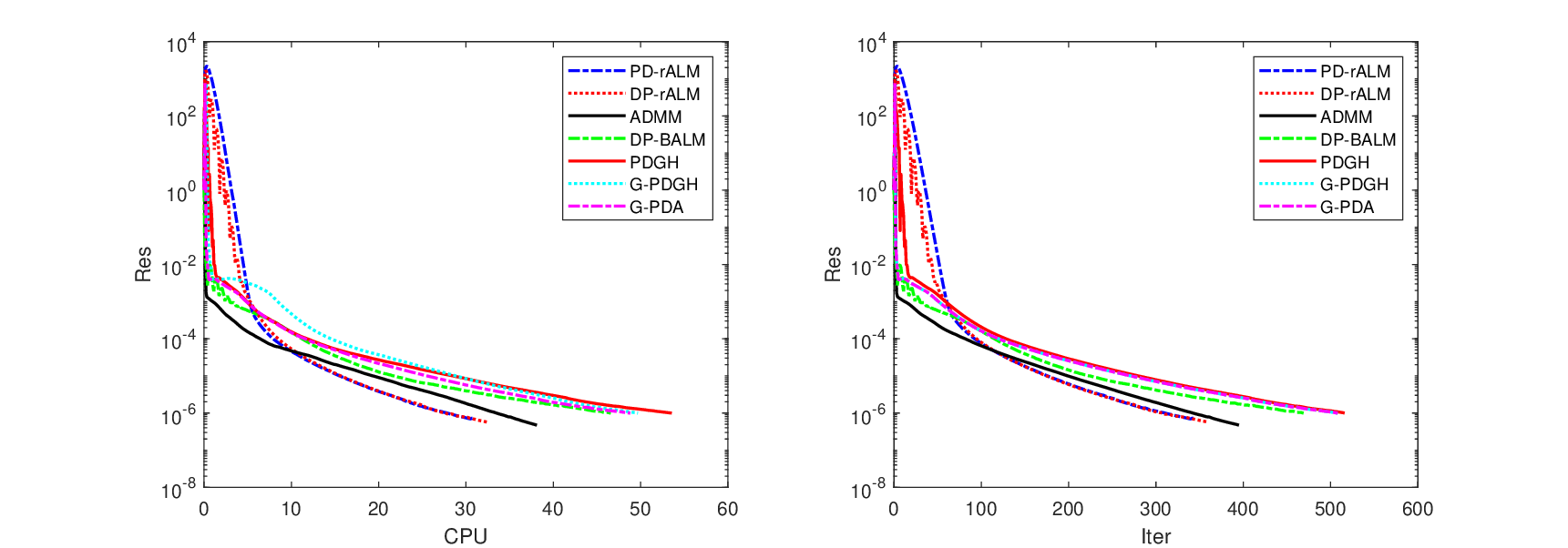}}\vspace{-0.4cm}
\caption{Comparative convergence curves of Res versus  CPU   and Iter,  respectively.} \label{Comp-curv}
   \end{minipage}
\end{figure}

{We report some comparative results of several state-of-the-art algorithms in Table \ref{TabProb-b} under different tolerances $(\epsilon_1,\epsilon_2)$.}
 The original columns of $D$, along with the low-rank and sparse components decomposed by different algorithms under $(\epsilon_1,\epsilon_2)=(10^{-5}, 10^{-6})$,
 are   shown in Figure  \ref{Fig4-2}, {and Figure \ref{Comp-curv} shows the comparative convergence curves of   Res(k) along with the increase  of  the CPU time  and the number of iterations.     From  Table \ref{TabProb-b} and Figure \ref{Comp-curv},  it can be seen that  both PD-rALM and DP-rALM perform better than other comparative methods in terms of the number of iterations and CPU time, and these two methods} can effectively fill  in
occluded regions of the image, corresponding to shadows. In the low-rank component
{$L$ as shown in Figure  \ref{Fig4-2}},   shadows  under different lighting conditions are removed and filled in with the most consistent low-rank features from
the eigenfaces.

\paragraph{Acknowledgements} {The authors would like to thank the editor and anonymous referees for their valuable comments and suggestions, which have significantly improved the quality of the paper.}

\paragraph{Funding}
This work  was supported by    Guangdong Basic and Applied Basic Research Foundation (2023A1515012405),  Shaanxi Fundamental Science Research Project for Mathematics and Physics (22JSQ001), and National Natural Science Foundation of China (Grants 12071398 {and 52372397}).

\paragraph{Data availability}
The link https://github.com/pzheng4218/new-ALM-in-ML provides the matlab codes for experiments.

\paragraph{Conflict of interest}
The authors declared that,  they do not have any commercial or associative interest that represents a conflict of interest in connection with  this paper.

\section{Appendix:   discussions on two new PDHG}
In this appendix, we   discuss two new types of PDHG algorithm without relaxation step   for  solving the following convex-concave saddle-point problem
\begin{equation} \label{Sec4-Prob}
\min\limits_{\m{x}\in\C{X}}\max\limits_{\m{y}\in\C{Y}}  \Phi(\m{x,y}):= \theta_1(\m{x})-\m{y}\tr A\m{x}-\theta_2(\m{y}),
\end{equation}
or, equivalently, the composite problem
$
\min\limits_{\m{x}\in\C{X}} \big\{\theta_1(\m{x}) + \theta_2^*(-A\m{x})\big\},
$
where $\C{X}\subseteq \C{R}^n, \C{Y}\subseteq \C{R}^m$ are closed convex sets,   both $\theta_1: {\C{X}}\rightarrow\C{R}$ and $ \theta_2: {\C{Y}}\rightarrow\C{R}$ are convex  but possibly nonsmooth functions, $\theta_2^*$ is the conjugate function of $\theta_2$, and $A\in\C{R}^{m\times n}$ is a given data. A lot of practical examples can be reformulated as   special cases of (\ref{Sec4-Prob}), see e.g. \cite[Section 5]{JWZC20}.  {Note that Problem (\ref{Sec4-Prob}) can  reduce to the dual of (\ref{Sec1-Prob1}) by letting $\theta_2=-\lambda\tr b,\m{y}=\lambda$ and $\C{Y}=\Lambda.$} Hence, the convergence results also hold for the previous   P-rALM.    Throughout the forthcoming discussions, the solution set of (\ref{Sec4-Prob}) is assumed to be nonempty.

The original PDHG proposed in \cite{Zhuchan08} is to solve some TV image restoration  problems. Extending it to the problem (\ref{Sec4-Prob}), we get the following scheme:
\[\left\{\begin{array}{l}
\m{x}^{k+1}=\arg \min\limits_{\m{x}\in\C{X}}  \Phi(\m{x},\m{y}^k)+\frac{r}{2}\|\m{x}-\m{x}^k\|^2, \\
\m{y}^{k+1}=\arg \max\limits_{\m{y}\in\C{Y}}  \Phi(\m{x}^{k+1},\m{y})-\frac{s}{2}\|\m{y}-\m{y}^k\|^2,
\end{array}\right.
\]
where $r,s$ are positive scalars. He, et al. \cite{hmy20} pointed out that   convergence of the above PDHG can be ensured if $\theta_1$ is strongly convex and $rs>\rho(A\tr A)$.   To weaken these convergence conditions, (e.g., the function $\theta_1$ is only convex and the parameters $r,s$ do not depend on  $\rho(A\tr A)$), we   develop a  novel PDHG (N-PDHG1) as   follows.
\begin{flushleft}
\centering\fbox{
	\parbox{0.89\textwidth}{
\[
   \begin{array}{lll}
 \hspace*{-.1in}\textrm{\bf Initialize } (\m{x}^0,\lambda^0) \textrm{ and choose }
  r>0,~ Q\succ\m{0}; \\
 \hspace*{-.1in} \textrm{\bf While } \textrm{the stopping criterion is not satisfied, }  \textrm{\bf do}\\
\hspace*{.1in}     \m{x}^{k+1}=\arg \min\limits_{\m{x}\in\C{X}}  \Phi(\m{x},\m{y}^k)  + \frac{1}{2}\left\|\m{x}-\m{x}^{k}\right\|^2_{rA\tr A+ Q};\\
\hspace*{.1in}  \m{y}^{k+1}=\arg \max\limits_{\m{y}\in\C{Y}}  \Phi(2\m{x}^{k+1}-\m{x}^k,\m{y})-\frac{1}{2r}\|\m{y}-\m{y}^k\|^2;\\
 \hspace*{-.1in} \textrm{\bf End while}
\end{array}
\]}}\end{flushleft}
Another  related  algorithm, called N-PDHG2, is just to modify the final subproblem of  PDHG, whose framework is  described in the next box.  A similar quadratic term was  adopted in \cite{HeDwang16} to solve a special case of   (\ref{Sec4-Prob}).   {We can observe that N-PDHG1 has certain connections  with   P-rALM, since the first-order optimality conditions  of their involved subproblems are reformulated as similar variational inequalities with the same block matrix $H$, see (\ref{b-3}) and the next (\ref{jb11-3}). Actually, their $\m{x}$-subproblems  enjoy the same proximal term. Another observation is that N-PDHG2 is developed from N-PDHG1 by just modifying the involved proximal parameters, and one of their subproblems could enjoy a proximity operator directly.}
\begin{flushleft}
\centering\fbox{
	\parbox{0.89\textwidth}{
\[
   \begin{array}{lll}
 \hspace*{-.1in}\textrm{\bf Initialize } (\m{x}^0,\lambda^0) \textrm{ and choose }
  r>0,~ Q\succ\m{0}; \\
 \hspace*{-.1in} \textrm{\bf While } \textrm{the stopping criterion is not satisfied, }  \textrm{\bf do}\\
\hspace*{.1in}     \m{x}^{k+1}=\arg \min\limits_{\m{x}\in\C{X}}  \Phi(\m{x},\m{y}^k)  + \frac{r}{2}\left\|\m{x}-\m{x}^{k}\right\|^2;\\
\hspace*{.1in}  \m{y}^{k+1}=\arg \max\limits_{\m{y}\in\C{Y}}  \Phi(2\m{x}^{k+1}-\m{x}^k,\m{y})-\frac{1}{2}\|\m{y}-\m{y}^k\|^2_{A A\tr/r+ Q};\\
 \hspace*{-.1in} \textrm{\bf End while}
\end{array}
\]}}\end{flushleft}

\subsection{Sublinear convergence under general  convex assumption}
Due to the close similarities between the two algorithms mentioned above, we will only analyze  the convergence properties of N-PDHG1 under general convex assumptions in the following section before proceeding over the second algorithm's convergence. For convenience, we denote $\C{U}:=\C{X}\times\C{Y}$ and
\[
\theta(\m{u})=\theta_1(\m{x})+\theta_2(\m{y}),\quad\m{u}=\left(\begin{array}{c}
 \m{x}  \\   \m{y}
\end{array}\right),\quad
\m{u}^k=\left(\begin{array}{c}
 \m{x}^k  \\   \m{y}^k
\end{array}\right)\quad\textrm{and} \quad M=\left[\begin{array}{ccccccc}
  \m{0} &&  -A\tr  \\
   A   &&    \m{0}
\end{array}\right].
\]

\begin{lemma} \label{jx1-optimal}
The sequence  $\{ \m{u}^k \}$   generated by  N-PDHG1 satisfies
\begin{equation}\label{jb11-3}
\m{u}^{k+1}\in \C{U},~ \theta(\m{u})- \theta(\m{u}^{k+1}) +\big\langle \m{u}-\m{u}^{k+1}, M\m{u}\big\rangle \geq   \big\langle\m{u}-\m{u}^{k+1}, H(\m{u}^k-\m{u}^{k+1}) \big\rangle
\end{equation}
for any $\m{u}\in \C{U}$, where $H$ is given by (\ref{H}).
 Moreover, we have
\begin{eqnarray}\label{jb11-4}
  \theta(\m{u})- \theta(\m{u}^{k+1}) +\big\langle \m{u}-\m{u}^{k+1}, M\m{u}\big\rangle
 \geq   \frac12\Big( \big\|\m{u}-\m{u}^{k+1}\big\|^2_H-\big\|\m{u}-\m{u}^k\big\|^2_H\Big)+\frac12 \big\|\m{u}^k-\m{u}^{k+1}\big\|^2_H.
\end{eqnarray}
\end{lemma}
 Proof.
According to the first-order optimality condition of the $\m{x}$-subproblem in N-PDHG1, we have $\m{x}^{k+1}\in\C{X}$ and
\begin{eqnarray}\label{jb1-5-0}
   \theta_1(\m{x})- \theta_1(\m{x}^{k+1})  +  \left\langle\m{x}-\m{x}^{k+1}, -A\tr \m{y}^k+\big(rA\tr A+Q\big) (\m{x}^{k+1}-\m{x}^k)\right\rangle\geq  0,   ~ \forall \m{x} \in\C{X},\end{eqnarray}
 that is,
\begin{eqnarray}\label{jb1-5}
&&\quad \theta_1(\m{x})- \theta_1(\m{x}^{k+1})  +  \left\langle\m{x}-\m{x}^{k+1}, -A\tr \m{y}^{k+1} \right\rangle\nonumber\\
&&\geq    \Big\langle\m{x}-\m{x}^{k+1},  \big(rA\tr A+Q\big) (\m{x}^k-\m{x}^{k+1})+ A\tr (\m{y}^k-\m{y}^{k+1})  \Big\rangle.
 \end{eqnarray}
Similarly, we have   $\m{y}^{k+1}\in \C{Y}$ and
\begin{eqnarray}\label{jb1-5-1}
   \theta_2(\m{y})- \theta_2(\m{y}^{k+1})  +  \Big\langle\m{y}-\m{y}^{k+1}, A  (2\m{x}^{k+1}-\m{x}^k)+\frac{1}{r} (\m{y}^{k+1}-\m{y}^k)\Big\rangle\geq  0,   ~ \forall \m{y} \in\C{Y}, \end{eqnarray}
 that is,
\begin{eqnarray}\label{jb1-6}
  \theta_2(\m{y})- \theta_2(\m{y}^{k+1})  +  \Big\langle\m{y}-\m{y}^{k+1}{,}   A\tr \m{x}^{k+1} \Big\rangle \geq    \Big\langle\m{y}-\m{y}^{k+1},  A (\m{x}^k-\m{x}^{k+1})+ \frac{1}{r} (\m{y}^k-\m{y}^{k+1})  \Big\rangle.
 \end{eqnarray}
 Combine the   inequalities (\ref{jb1-5})-(\ref{jb1-6}) and  the structure of $H$ given by (\ref{H}) to have
 \[
  \theta(\m{u})- \theta(\m{u}^{k+1}) +\Big\langle \m{u}-\m{u}^{k+1}, M\m{u}^{k+1}\Big\rangle \geq   \Big\langle\m{u}-\m{u}^{k+1}, H(\m{u}^k-\m{u}^{k+1}) \Big\rangle,
 \]
 which together with the the property $\left\langle \m{u}-\m{u}^{k+1}, M(  \m{u}-\m{u}^{k+1}) \right\rangle=0$ confirms (\ref{jb11-3}). Then, the  inequality (\ref{jb11-4}) is obtained {by applying (\ref{jb11-3}) and    the identity in (\ref{identabcd}).}
 $\hfill \blacksquare$

 Now, we    discuss  the global convergence and sublinear convergence rate  of N-PDHG1. Let $\m{u}^*=(\m{x}^*;\m{y}^*)\in\C{U}$ be a solution point of the problem (\ref{Sec4-Prob}). Then, it holds
 \[
  {\Phi(\m{x}^*,\m{y})\leq \Phi(\m{x}^*,\m{y}^*)\leq \Phi(\m{x},\m{y}^*), \quad\forall \m{x}\in\C{X}, \m{y}\in\C{Y}},
 \]
 namely,
 \[
   \left \{\begin{array}{lllll}
\m{x}^*\in\C{X}, &\theta_1(\m{x})- \theta_1(\m{x}^*) &+ &\langle\m{x}-\m{x}^*,  -A\tr \m{y}^*\rangle\geq  0,  &\forall \m{x} \in\C{X},\\
\m{y}^*\in\C{Y}, &\theta_2(\m{y})- \theta_2(\m{y}^*) &+ &\langle\m{y}-\m{y}^*,   A \m{x}^*\rangle\geq  0,  &\forall \m{x} \in\C{Y}.
\end{array}\right.
\]
So, finding a solution point of  (\ref{Sec4-Prob}) amounts to finding $\m{u}^*\in \C{U}$ such that
\begin{equation}\label{jb11-7}
\m{u}^*\in \C{U},~~~ \theta(\m{u})- \theta(\m{u}^*) +\left\langle \m{u}-\m{u}^*, M\m{u}^*\right\rangle \geq  0, \quad\forall \m{u}\in \C{U}.
\end{equation}
Setting $\m{u}:=\m{u}^*$ in (\ref{jb11-4}) together with (\ref{jb11-7}) gives
\begin{equation}\label{jb11-key}
\big \|\m{u}^*-\m{u}^{k+1}\big\|^2_H\leq\big\|\m{u}^*-\m{u}^k\big\|^2_H -
     \big\|\m{u}^k-\m{u}^{k+1}\big\|^2_H,
\end{equation}
that is, the sequence generated by N-PDHG1 is contractive, and thus N-PDHG1 converges globally.  The last  inequality together with the analysis  of P-rALM indicates that N-PDHG1 with a relaxation step also converges, and   the sublinear convergence rate of N-PDHG1 is similar to the proof of P-rALM. Note that the convergence of   N-PDHG1 does not need the strong  convexity of $\theta_1$ and allows more flexibility for choosing the proximal parameter $r$.

Finally, it is not difficult  from the first-order optimality conditions of the involved subproblems in N-PDHG2 that
\[
\m{u}^{k+1}\in \C{U},~~ \theta(\m{u})- \theta(\m{u}^{k+1}) +\Big\langle \m{u}-\m{u}^{k+1}, M\m{u}\Big\rangle \geq   \Big\langle\m{u}-\m{u}^{k+1}, \widetilde{H}(\m{u}^k-\m{u}^{k+1}) \Big\rangle
\]
for any $\m{u}\in \C{U}$, where
\[
\widetilde{H}=\left[\begin{array}{ccccccc}
 r\m{I}  &&   A\tr  \\
   A   &&    \frac{1}{r}AA\tr+Q
\end{array}\right]
\]
and   $\widetilde{H}$ is positive definite for any $r>0$ and $Q\succ\m{0}$. So, N-PDHG2 also converges globally with a sublinear convergence rate. This matrix $\widetilde{H}$ is what we discussed in Section 1 and could reduce to that in \cite{Heyuan21} with $Q=\delta\m{I}$ for any $\delta>0$.

\subsection{Linear convergence under strongly convexity assumption}
The linear convergence rate of N-PDHG1 will be investigated in this subsection     under the following   assumptions:
\begin{itemize}
  \item [(a1)]
  The matrix $A$ has full row rank and $\C{X}=\C{R}^n;$
  \item [(a2)]
  The function $\theta_1$ is strongly convex with modulus $\nu>0$ and $\nabla\theta_1$ is Lipschitz continuous with constant $L_{\theta_1}>0$.
\end{itemize}

From the second part  of   (a2) and  the first-order optimality condition of $\m{x}^{k+1}$-subproblem in N-PDHG1, we have
\begin{equation}\label{sec62-1}
-\nabla\theta_1(\m{x}^{k+1})=-A\tr \m{y}^k+\big(rA\tr A+Q\big) (\m{x}^{k+1}-\m{x}^k).
\end{equation}
Together with this equation and the first part  of (a2), it holds
\begin{eqnarray*}
&&\theta_1(\m{x})-\theta_1(\m{x}^{k+1})\geq \big\langle\m{x}-\m{x}^{k+1}, \nabla\theta_1(\m{x}^{k+1})\big\rangle +\frac{\nu}{2}\big\|\m{x}-\m{x}^{k+1} \big\|^2\Rightarrow\nonumber\\
&&\theta_1(\m{x})-\theta_1(\m{x}^{k+1}) +\big\langle\m{x}-\m{x}^{k+1}, -A\tr \m{y}^{k+1}\big\rangle \geq \frac{\nu}{2}\big\|\m{x}-\m{x}^{k+1} \big\|^2+\nonumber\\
&&\qquad \big\langle\m{x}-\m{x}^{k+1},  \big(rA\tr A+Q\big) (\m{x}^k-\m{x}^{k+1})+ A\tr (\m{y}^k-\m{y}^{k+1})  \big\rangle,
\end{eqnarray*}
which implies that   $\frac{\nu}{2}\left\|\m{x}-\m{x}^{k+1} \right\|^2$ will be added to the right-hand-side of (\ref{jb11-3})  and finally
\begin{equation}\label{sec62-2}
\big \|\m{u}^*-\m{u}^{k+1}\big\|^2_H\leq\big\|\m{u}^*-\m{u}^k\big\|^2_H -
     \big\|\m{u}^k-\m{u}^{k+1}\big\|^2_H - \nu\big\|\m{x}^*-\m{x}^{k+1} \big\|^2.
\end{equation}
Note that the equation (\ref{sec62-1}) can be equivalently  rewritten as
\begin{equation}\label{sec62-3}
A\tr \m{y}^{k+1}= \nabla\theta_1(\m{x}^{k+1})+A\tr (\m{y}^{k+1}-\m{y}^k)+\big(rA\tr A+Q\big) (\m{x}^{k+1}-\m{x}^k).
\end{equation}
Besides, the solution $(\m{x}^*;\m{y}^*)$   satisfies
\begin{equation}\label{sec62-4}
\nabla\theta_1(\m{x}^*)= A\tr \m{y}^*.
\end{equation}
Combining the equations (\ref{sec62-3})-(\ref{sec62-4}) together with (a1)-(a2) is to obtain
\begin{eqnarray*}
  &&\sigma_A\big\|\m{y}^{k+1}-\m{y}^*\big\|^2 \leq  \big\|A\tr(\m{y}^{k+1}-\m{y}^*)\big\|^2 \\
&=& \big\|\nabla\theta_1(\m{x}^{k+1})- \nabla\theta_1(\m{x}^*) + A\tr (\m{y}^{k+1}-\m{y}^k)+\big(rA\tr A+Q\big) (\m{x}^{k+1}-\m{x}^k)\big\|^2 \\
&\leq & 3\Big\{\big\|\nabla\theta_1(\m{x}^{k+1})- \nabla\theta_1(\m{x}^*) \big\|^2 +\big\|  A\tr (\m{y}^{k+1}-\m{y}^k) \big\|^2+\big\| \big(rA\tr A+Q\big) (\m{x}^{k+1}-\m{x}^k)\big\|^2\Big\}\\
&\leq & 3\Big\{L_{\theta_1}^2\big\| \m{x}^{k+1} -  \m{x}^*  \big\|^2 +\|A \|^2\big\|   \m{y}^{k+1}-\m{y}^k  \big\|^2+\big\| \big(rA\tr A+Q\big) (\m{x}^{k+1}-\m{x}^k)\big\|^2\Big\},
\end{eqnarray*}
where $\sigma_A>0$ denotes the smallest   eigenvalue of $AA\tr$ due to (a1). So, we have
\begin{eqnarray}
    \big \|\m{u}^*-\m{u}^{k+1}\big\|^2_H&\leq& \|H\| \Big\{\big \|\m{x}^*-\m{x}^{k+1}\big\|^2+\big \|\m{y}^*-\m{y}^{k+1}\big\|^2\Big\}\nonumber\\
&\leq&  \|H\| \Big\{(1+3L_{\theta_1}^2\sigma_A^{-1})\big \|\m{x}^*-\m{x}^{k+1}\big\|^2+3\sigma_A^{-1}\|A \|^2\big\|   \m{y}^{k+1}-\m{y}^k  \big\|^2 \nonumber\\
&  & \qquad~  + 3\sigma_A^{-1}\big\| \big(rA\tr A+Q\big) (\m{x}^{k+1}-\m{x}^k)\big\|^2\Big\}.\label{sec62-5}
\end{eqnarray}
By the structure of $H$ and the Young's inequality,  it holds that
\begin{eqnarray}
 &&   \big \|\m{u}^k-\m{u}^{k+1}\big\|^2_H\nonumber\\
&=& \big\| \big(rA\tr A+Q\big) (\m{x}^{k+1}-\m{x}^k)\big\|^2+ \frac1r \big\|   \m{y}^{k+1}-\m{y}^k  \big\|^2+2\big\langle \m{x}^{k+1}-\m{x}^k, A\tr(\m{y}^{k+1}-\m{y}^k)\big\rangle \nonumber\\
& \geq &   \big\| \big(rA\tr A+Q\big) (\m{x}^{k+1}-\m{x}^k)\big\|^2+ \frac1r \big\|   \m{y}^{k+1}-\m{y}^k  \big\|^2\nonumber\\
&&-\Big\{ \delta_0\big\|  \m{x}^{k+1}-\m{x}^k \big\|^2+\frac{1}{\delta_0}\|A\tr A\|\big\|   \m{y}^{k+1}-\m{y}^k   \big\|^2\Big\}, \nonumber\\
& \geq &\big\| \big(rA\tr A+Q\big) (\m{x}^{k+1}-\m{x}^k)\big\|^2- \delta_0\big\|  \m{x}^{k+1}-\m{x}^k \big\|^2+\Big( \frac1r-\frac{ \|A\tr A\| }{\delta_0}\Big)\big\|   \m{y}^{k+1}-\m{y}^k   \big\|^2,\quad\label{sec62-6}
\end{eqnarray}
where $ \delta_0 \in( r\|A\tr A\|, \|rA\tr A+Q\|^2)$ exists for   proper choices of $r$ and $Q$. Now, let
\[
\delta^k= \min\left\{ \frac{\nu}{(1+3L_{\theta_1}^2\sigma_A^{-1})\|H\|},~
\frac{\delta_0-r\|A\tr A\|}{3r\delta_0\sigma_A^{-1}\|A\|^2\|H\|},~
\frac{\|rA\tr A+Q\|^2-\delta_0}{3 \sigma_A^{-1}\|H\|\left\| \left(rA\tr A+Q\right) (\m{x}^{k+1}-\m{x}^k)\right\|^2}
\right\}.
\]
Then, combining the above inequalities (\ref{sec62-2}) and (\ref{sec62-5})-(\ref{sec62-6}), we can deduce
\begin{eqnarray}
 &&  (1+\delta^k)\big \|\m{u}^*-\m{u}^{k+1}\big\|^2_H-\big\|\m{u}^*-\m{u}^k\big\|^2_H \nonumber\\
& \leq&  \delta^k\big\|\m{u}^*-\m{u}^{k+1}\big\|^2_H-
     \big\|\m{u}^k-\m{u}^{k+1}\big\|^2_H - \nu\big\|\m{x}^*-\m{x}^{k+1} \big\|^2\nonumber\\
&\leq& \Big\{\delta^k(1+3L_{\theta_1}^2 \sigma_A^{-1})\|H\|-\nu \Big\}\big\|\m{x}^*-\m{x}^{k+1} \big\|^2+\Big\{3\delta^k\sigma_A^{-1}\|A \|^2\|H\|-\frac1r+\frac{ \|A\tr A\|}{\delta_0}\Big\}\big\|   \m{y}^{k+1}-\m{y}^k   \big\|^2\nonumber\\
&  & +\big(3\delta^k\sigma_A^{-1}\|H\|-1\big)\big\| \big(rA\tr A+Q\big) (\m{x}^{k+1}-\m{x}^k)\big\|^2+\delta_0\big\|  \m{x}^{k+1}-\m{x}^k \big\|^2.\label{sec62-7}
\end{eqnarray}
Observing from the definition of $\delta^k$,  it holds
\[
   \left \{\begin{array}{lllll}
\delta^k(1+3L_{\theta_1}^2\sigma_A^{-1})\|H\|-\nu\leq 0,\\
 3\delta^k\sigma_A^{-1}\|A \|^2\|H\|-\frac1r+\frac{\|A\tr A\|}{\delta_0}\leq 0,\\
 \left(3\delta^k\sigma_A^{-1}\|H\|-1\right)\left\| \left(rA\tr A+Q\right) (\m{x}^{k+1}-\m{x}^k)\right\|^2+\delta_0\left\|  \m{x}^{k+1}-\m{x}^k \right\|^2\leq 0,
\end{array}\right.
\]
and finally ensures the following $Q$-linear convergence rate:
\[
\big\|\m{u}^*-\m{u}^{k+1}\big\|^2_H\leq
\frac{1}{ 1+\delta^k }\big\|\m{u}^*-\m{u}^k\big\|^2_H.
\]
The above analysis also indicates that our proposed P-rALM for solving the problem (\ref{Sec1-Prob1})  will converge  $Q$-linearly under the similar assumptions that $\theta$ is strongly convex, its gradient $\nabla\theta$ is Lipschitz continuous, the matrix $A$ has full row rank and $\C{X}=\C{R}^n$.

\subsection{Linear convergence under the error bound condition}
In this section, we use $\partial f(x)$ to denote the sub-differential of the  convex function $f$ at $x$.   $f$ is said to be a piecewise linear multifunction if its graph $Gr(f):=\{(x,y)\mid y\in f(x)\}$  is a union of finitely many polyhedra.  The    projection operator  $\C{P}_{\C{C}}(x)$   is   nonexpansive,  i.e.,
\begin{equation}\label{Apend-0}
\|\C{P}_{\C{C}}(x)-\C{P}_{\C{C}}(z)\|\leq \|x-z\|,\quad\forall x,z\in\C{R}^n.
\end{equation}
Given $H\succ \m{0}$, we define $\operatorname{dist}_{H}(x, \C{C}):=\min\limits_{z\in\C{C}}\|x-z\|_H$. When $H=\m{I}$, we simply denote it by  $\operatorname{dist}(x, \C{C})$.
For any $\m{u}\in\C{U}$ and $\alpha>0$, we  define
 \begin{equation}\label{Apend-1}
e_{\C{U}}(\m{u},\alpha):=\left(\begin{array}{c}
 e_{\C{X}}(\m{u},\alpha):= \m{x}-\C{P}_{\C{X}}\left[\m{x} - \alpha(\xi_{\m{x}} -A\tr \m{y}{)}\right]\\
e_{\C{Y}}(\m{u},\alpha):= \m{y}-\C{P}_{\C{Y}}\left[\m{y} - \alpha(\xi_{\m{y}} +A  \m{x}{)}\right]
\end{array}\right),
\end{equation}
where $\xi_{\m{x}}\in\partial\theta_1(\m{x}),\xi_{\m{y}}\in\partial\theta_2(\m{y}).$
 Note that a point
 \[
 \m{u}^*\in \C{U}^* = \big\{\hat{\m{u}}\in  \C{U}\mid \operatorname{dist}\left(\mathbf{0},e_{\C{U}}(\hat{\m{u}},\alpha)\right)
 =0\big\}
 \]
 is the solution of (\ref{Sec4-Prob}) if and only if $ e_{\C{U}}(\m{u}^*,\alpha)=\m{0}$.    Different from the assumptions {(a1)-(a2)},   we next investigate the   linear convergence  rate of N-PDHG1 under    an error bound condition in terms of the mapping $e_{\C{U}}(\m{u},1)$:
\begin{itemize}
  \item [(a3)]
  Assume that there exists a constant $\zeta>0$ such that
\begin{equation}\label{Apend-2}
\operatorname{dist}\left(\m{u},\C{U}^*\right)\leq \zeta \operatorname{dist}\left(\mathbf{0},e_{\C{U}}(\m{u},1)\right), ~~\forall \m{u}\in\C{U}.
\end{equation}
\end{itemize}

The condition  (\ref{Apend-2}) is  generally weaker than the strong convexity assumption
and hence can be satisfied by some problems that have non-strongly convex objective functions.  Note that if   the sub-differentials $\partial\theta_1(\m{x})$ and $\partial\theta_2(\m{y})$ are piecewise linear multifunctions and
the constraint sets $\C{X}, \C{Y}$ are polyhedral, then
  both $\C{P}_{\C{X}}$ and $\C{P}_{\C{Y}}$ are  piecewise linear multifunctions by \cite[Prop. 4.1.4]{Fpa03} and hence $e_{\C{U}}(\m{u},\alpha)$ is also a piecewise linear multifunction.   Followed by Robinson's continuity property \cite{Ro81} for polyhedral multifunctions, the assumption (a3) holds   automatically.
For convenience of the sequel analysis, we denote
 \begin{equation}\label{Apend-3}
 \C{Q} =\left[\begin{array}{ccccccc}
   (rA\tr A+Q)\tr(rA\tr A+Q)+ A\tr A &&  \m{0}   \\
      \m{0} &&    \frac1r\m{I}+AA\tr
\end{array}\right].
\end{equation}
It is easy to check that $\C{Q}$ is symmetric positive definite because $\|\m{u}\|^2_{\C{Q}}>0$ for any $\m{u}\neq \m{0}$. By   equivalent expressions for the first-order optimality conditions (\ref{jb1-5-0}) and (\ref{jb1-5-1}) together with the structure of $\C{Q}$,  we have the following estimation on the distance {between  $\m{0}$ and}   $e_{\C{U}}(\m{u}^{k+1},1)$, which follows   the similar proof as that in \cite[Sec. 2.2]{BCLXu21}.
\begin{lemma}\label{Apend-4}
Let  $\C{Q}$ be given in (\ref{Apend-3}). Then, the iterates  generated by N-PDHG1 satisfy
 \begin{equation}\label{Apend-5}
\operatorname{dist}^2\big(\mathbf{0},e_{\C{U}}(\m{u}^{k+1},1)\big)\leq  2 \big\|\m{u}^k-\m{u}^{k+1}\big\|_{\C{Q}}^2.
\end{equation}
\end{lemma}
 Proof. The first-order optimality condition in (\ref{jb1-5-0}) implies
 \[
\m{x}^{k+1}=\C{P}_{\C{X}}\left\{\m{x}^{k+1}-\Big[\xi_{\m{x}}^{k+1} -A\tr \m{y}^k+\big(rA\tr A+Q\big) (\m{x}^{k+1}-\m{x}^k)\Big]\right\}.
 \]
Combine it    with the definition of $\operatorname{dist}_H(\cdot,\cdot)$ and the property in (\ref{Apend-0}) to obtain
 \begin{eqnarray}
 && \operatorname{dist}^2\left(\mathbf{0},e_{\C{X}}(\m{u}^{k+1},1)\right) = \operatorname{dist}^2\left(\m{x}^{k+1}, \C{P}_{\C{X}}\left\{\m{x}^{k+1}-\big[\xi_{\m{x}}^{k+1} -A\tr \m{y}^{k+1}\big]\right\}\right) \nonumber\\
& \leq& \Big\|A\tr(\m{y}^k-\m{y}^{k+1}) + \big(rA\tr A+Q\big) (\m{x}^k-\m{x}^{k+1})\Big\|^2  \nonumber\\
&\leq&  2\Big( \big\|A\tr(\m{y}^k-\m{y}^{k+1})\big\|^2+ \big\|\big(rA\tr A+Q\big) (\m{x}^k-\m{x}^{k+1})\big\|^2\Big)=2 \big\| \m{u}^k-\m{u}^{k+1}\big\|^2_{\C{Q}_1}, \label{Apend-6}
\end{eqnarray}
where $\C{Q}_1=\operatorname{diag}\left((rA\tr A+Q)\tr(rA\tr A+Q),AA\tr\right).$  Similarly, we have from (\ref{jb1-5-1}) that
 \[
\m{y}^{k+1}=\C{P}_{\C{Y}}\left\{\m{y}^{k+1}-\Big[\xi_{\m{y}}^{k+1}+ A  (2\m{x}^{k+1}-\m{x}^k)+\frac{1}{r} (\m{y}^{k+1}-\m{y}^k)\Big]\right\}
 \]
and
 \begin{eqnarray}
 && \operatorname{dist}^2\big(\mathbf{0},e_{\C{Y}}(\m{u}^{k+1},1)\big) = \operatorname{dist}^2\Big(\m{y}^{k+1}, \C{P}_{\C{Y}}\big\{\m{y}^{k+1}-\big[\xi_{\m{y}}^{k+1} +A  \m{x}^{k+1}\big]\big\}\Big) \nonumber\\
& \leq& \Big\|A (\m{x}^k-\m{x}^{k+1}) +\frac{1}{r} (\m{y}^k-\m{y}^{k+1})\Big\|^2  \nonumber\\
&\leq&  2\Big( \big\|A (\m{x}^k-\m{x}^{k+1})\big\|^2+ \big\|\frac{1}{r} (\m{y}^k-\m{y}^{k+1})\big\|^2\Big)=2 \big\| \m{u}^k-\m{u}^{k+1}\big\|^2_{\C{Q}_2},\label{Apend-7}
\end{eqnarray}
where $\C{Q}_2=\operatorname{diag}\left( A\tr A,\frac1r\m{I}\right).$
The   inequalities (\ref{Apend-6})-(\ref{Apend-7}) immediately  ensure  (\ref{Apend-5}) due to the relation $\C{Q}=\C{Q}_1+\C{Q}_2.$
 $\hfill \blacksquare$

 Based on Lemma \ref{Apend-4} and the conclusion (\ref{jb11-key}), we next provide   a  global linear convergence rate of N-PDHG1 with the aid of the notations   $\lambda_{\min}(H)$ and $ \lambda_{\max}(H)$ which  denote  the smallest  and largest eigenvalue of the positive definite matrix $H$,  respectively.
\begin{theorem}\label{Apend-8}
Let  $\C{Q}$ be given in (\ref{Apend-3}).  Then,
there exists a constant $\zeta>0$ such that  the   iterates  generated by N-PDHG1  {satisfy}
\begin{equation}\label{Apend-9}
\operatorname{dist}^2_{H}(\m{u}^{k+1}, \C{U}^*)\leq \frac{1}{1+\hat{\zeta}  }\operatorname{dist}^2_{H}(\m{u}^k, \C{U}^*),
\end{equation}
where  the constant
$
\hat{\zeta}= \frac{\lambda_{\min}(H)}{2\zeta^2\lambda_{\max}(\C{Q})\lambda_{\max}(H)}>0.
$
\end{theorem}
\noindent{\bf Proof }
Because $ \C{U}^*$ is a closed convex  set, there exists a $\m{u}^*_k\in \C{U}^*$ satisfying
\begin{equation} \label{Apend-010}
\operatorname{dist}_{H}(\m{u}^k, \C{U}^*)
=\big\|\m{u}^k-\m{u}^*_k\big\|_{H}.
\end{equation}
 By  the condition (\ref{Apend-2}) and Lemma \ref{Apend-4}  there exists a constant $\zeta>0$ such that
\begin{equation} \label{Apend-10}
\operatorname{dist}^2\big(\m{u}^{k+1}, \C{U}^*\big)\leq 2\zeta^2\big\|\m{u}^k-\m{u}^{k+1}\big\|_{\C{Q}}^2\leq
\frac{2\zeta^2\lambda_{\max}(\C{Q})}{\lambda_{\min}(H)}\big\|\m{u}^k
-\m{u}^{k+1}\big\|_{H}^2.
\end{equation}
By the definition of $\operatorname{dist}_{H}(\cdot,\cdot)$, {we have}
\begin{equation} \label{Apend-11}
\frac{1}{\lambda_{\max}(H)}\operatorname{dist}^2_H\big(\m{u}^{k+1}, \C{U}^*\big)\leq \operatorname{dist}^2\big(\m{u}^{k+1}, \C{U}^*\big).
\end{equation}
Combine (\ref{Apend-10})-(\ref{Apend-11}) and (\ref{jb11-key}) to have
\begin{eqnarray*} \label{Sec33-1090}
&& \operatorname{dist}^2_{H}(\m{u}^{k+1}, \C{U}^*)\leq \big\|\m{u}^{k+1}-\m{u}^*_k\big\|_H^2 \nonumber\\
&\leq&  \big\|\m{u}^k-\m{u}^*_k\big\|_H^2 -  \big\|\m{u}^k-\m{u}^{k+1}\big\|_H^2\\
&\leq&{\operatorname{dist}^2_{H}}(\m{u}^k, \C{U}^*)-\frac{\lambda_{\min}(H)}{2\zeta^2\lambda_{\max}(\C{Q})\lambda_{\max}(H)} \operatorname{dist}^2_H\big(\m{u}^{k+1}, \C{U}^*\big).
\end{eqnarray*}
Rearranging the above inequality is to confirm (\ref{Apend-9}).  $\hfill    \blacksquare$

\begin{corollary} \label{13131}
  Let $\hat{\zeta}>0$ be given  in Theorem \ref{Apend-8} and the  sequence   $\{\m{u}^k\}$ be generated by
N-PDHG1. Then, there exists a point $\m{u}^\infty \in \C{U}^*$
such that
\begin{equation}\label{v-lin-conv}
\big\|\m{u}^k-\m{u}^\infty \big\|_{H}\leq C  \epsilon^k,
\end{equation}
where
\[
C=\frac{2\operatorname{dist}_{H}(\m{u}^0, \C{U}^*)}{1-\epsilon}>0\quad \mbox{and}\quad \epsilon=\frac{1}{\sqrt{1+\hat{\zeta}}}\in(0,1).
\]
\end{corollary}
\noindent{\bf Proof }
Let ${\m{u}^*}\in \C{U}^*$ such that (\ref{Apend-010}) holds and let
\begin{equation}\label{abbc}
\m{u}^{k+1}=\m{u}^{k}+\m{d}^k.
\end{equation}
Then,   it follows from (\ref{jb11-key})  that
$
\left\|\m{u}^{k+1}-{\m{u}^*}\right\|_{H}\leq \left\|\m{u}^k-{\m{u}^*}\right\|_{H}
$
 which further implies
\begin{eqnarray}\label{d-linear}
\big\|\m{d}^k\big\|_H&=&\big\|\m{u}^{k+1}-\m{u}^k\big\|_{H}  \leq   \big\|\m{u}^{k+1}-{\m{u}^*}\big\|_{H}+
\big\|\m{u}^k-{\m{u}^*}\big\|_{H} \nonumber\\
&\leq & 2\big\|\m{u}^k-{\m{u}^*}\big\|_{H}
= 2\operatorname{dist}_{H}(\m{u}^k, \C{U}^*)\nonumber\\
&\leq &2 \epsilon^k \operatorname{dist}_{H}\big(\m{u}^0, \C{U}^*\big),
\end{eqnarray}
where the final inequality follows  from (\ref{Apend-9}).
Because  the sequence  $\{\m{u}^k\} $
generated by N-PDHG1
converges to a $ \m{u}^\infty  \in \C{U}^*$,
 we have from (\ref{abbc}) that $ \m{u}^\infty =\m{u}^k+\sum_{j=k}^{\infty}\m{d}^j $,
which by (\ref{d-linear}) indicates
\begin{eqnarray*}
\big\|\m{u}^k-\m{u}^\infty\big\|_{H}&\leq&\sum\limits_{j=k}^{\infty}\|\m{d}^j\|_{H}
\leq  2\operatorname{dist}_{H}(\m{u}^0, \C{U}^*)\sum\limits_{j=k}^{\infty}\epsilon^{j}  \\
&= &  2\operatorname{dist}_{H}(\m{u}^0, \C{U}^*) \epsilon^k \sum\limits_{j=0}^{\infty}\epsilon^{j} \leq   \epsilon^k \Big[2\operatorname{dist}_{H}(\m{u}^0, \C{U}^*)\frac{1}{1-\epsilon}\Big].
\end{eqnarray*}
So, the inequality (\ref{v-lin-conv}) holds, that is, $\m{u}^k$ converges $\m{u}^\infty$ R-linearly.  $\hfill  \blacksquare$

\begin{remark}
Consider the following general  saddle-point problem
\[
\min\limits_{\m{x}\in\C{X}}\max\limits_{\m{y}\in\C{Y}}  \Phi(\m{x,y}):= f(\m{x})+\theta_1(\m{x})-\m{y}\tr A\m{x}-\theta_2(\m{y}),
\]
or, equivalently, the composite problem
$
\min\limits_{\m{x}\in\C{X}} \big\{f(\m{x}) + \theta_1(\m{x}) + \theta_2^*(-A\m{x})\big\},$
where   $f(\m{x}): {\C{X}}\rightarrow\C{R}$ is a smooth convex function and its gradient is Lipschitz continuous
with constant $L_f$, and the remaining notations   have the same meanings as before. For this problem,   similar to the previous case 2 in Section \ref{twocases} we can develop the following iterative scheme
\[
   \left \{\begin{array}{lllll}
\m{x}^{k+1}=\arg \min\limits_{\m{x}\in\C{X}} \theta_1(\m{x})   +
 \big\langle \nabla f(\m{x}^k)-A\tr\m{y}^k, \m{x}\big\rangle+ \frac{1}{2}\left\|\m{x}-\m{x}^{k}\right\|^2_{rA\tr A+ Q},\\
\m{y}^{k+1}=\arg \max\limits_{\m{y}\in\C{Y}}  \Phi(2\m{x}^{k+1}-\m{x}^k,\m{y})-\frac{1}{2r}\|\m{y}-\m{y}^k\|^2.
\end{array}\right.
\]
Its global convergence and   linear  convergence rate can be also established   by   the above  analysis.
\end{remark}

\end{document}